\documentclass[a4paper,10pt]{scrartcl}
\usepackage[english]{babel}
\usepackage[utf8]{inputenc}
\usepackage{lmodern}
\usepackage{enumerate}
\usepackage{geometry}
\usepackage{amsfonts,amsmath,amssymb,amsopn}
\usepackage{amsthm}
\usepackage{mathtools}
\usepackage{stmaryrd}
\usepackage{nicefrac}
\usepackage{exscale}

\usepackage[numbers,square]{natbib}
\usepackage{url} 
\usepackage[colorlinks=true,pdfpagelabels,unicode]{hyperref}
\hypersetup{urlcolor=red, citecolor=blue} 

\allowdisplaybreaks 

\theoremstyle{plain}
\newtheorem{thmx}{Theorem}

\newtheoremstyle{theo}
	{3pt} 
	{3pt} 
	{\itshape} 
	{} 
		{\bfseries} 
	{\\} 
	{ } 
	{\thmname{#1}\thmnumber{ #2.}\thmnote{ - #3}} 
\theoremstyle{theo}

\newtheorem{definition}{Definition}[section]
\newtheorem{lemma}[definition]{Lemma}
\newtheorem{theorem}[definition]{Theorem}

\newtheorem{remark}[definition]{Remark}
\newtheorem{proposition}[definition]{Proposition}

\newenvironment{bew}{\begin{proof}[\bfseries Proof:]}{\end{proof}}

\DeclareMathOperator{\bomega}{\bar{\Omega}}
\DeclareMathOperator{\romega}{\partial\Omega}
 
\DeclareMathOperator{\D}{D}
\DeclareMathOperator{\intd}{d\!}
\DeclareMathOperator{\dive}{\nabla\cdot}
\DeclareMathOperator{\wto}{\rightharpoonup}
\DeclareMathOperator{\wsto}{\stackrel{\star}{\wto}}
\newcommand{\eps}{\varepsilon}

\newcommand{\divedot}{\cdot\!}
\newcommand{\F}{F_{\mu}}
\newcommand{\ms}{m_\star}
\newcommand{\mss}{m_{\star\star}}
\newcommand{\ts}{t_\star}
\newcommand{\tss}{t_{\star\star}}
\newcommand{\Tm}{T_{max}}
\newcommand{\GNI}{Gagliardo--Nirenberg inequality}

\newcommand{\intot}{\int_0^{t}\!}
\newcommand{\intomega}{\int_{\Omega}\!} 
\newcommand{\intinfomega}{\int_0^\infty\!\!\intomega}
\newcommand{\intromega}{\int_{\romega}\!} 

\newcommand{\Lo}[1][1]{L^{#1}(\Omega)} 
\newcommand{\W}[1][1,2]{W^{#1}(\Omega)}
\newcommand{\Co}[1][0]{C^{#1}(\Omega)}
\newcommand{\LSp}[2]{L^{#1\;\!}\!\left(#2\right)} 
\newcommand{\LSploc}[2]{L_{loc}^{#1}\!\left(#2\right)} 
\newcommand{\WSp}[2]{W^{#1}\!\left(#2\right)}
\newcommand{\CSp}[2]{C^{#1}\!\left(#2\right)}
\newcommand{\CSph}[2]{C^{#1}\!\,\left(#2\right)}
\newcommand{\CSploc}[2]{C_{loc}^{#1}\!\left(#2\right)}
\newcommand{\DA}[1][\alpha]{D\!\left(A^{#1}\right)} 
\newcommand{\DAr}[2]{D\!\left(A^{#1}_{#2}\right)}
\newcommand{\R}{\mathbb{R}}
\newcommand{\N}{\mathbb{N}}
\newcommand{\HP}{\mathcal{P}}
\newcommand{\HPr}{\mathcal{P}_r}
\newcommand{\nfrac}[2]{{\nicefrac{#1}{#2}}}

\author{Tobias Black\thanks{Institut f\"ur Mathematik, Universit\"at Paderborn, Warburger Str. 100, 33098 Paderborn, Germany; Email address: \mbox{tblack@math.upb.de}}}
\title{Eventual smoothness of generalized solutions to a singular chemotaxis-Stokes system}
\begin{document}
\maketitle
\begin{abstract}
\noindent
{\textbf{Abstract:} 
We study the chemotaxis-fluid system
\begin{align*}
\left\{
\begin{array}{r@{\,}c@{\,}c@{\ }l@{\quad}l@{\quad}l@{\,}c}
n_{t}&+&u\cdot\!\nabla n&=\Delta n-\nabla\!\cdot(\frac{n}{c}\nabla c),\ &x\in\Omega,& t>0,\\
c_{t}&+&u\cdot\!\nabla c&=\Delta c-nc,\ &x\in\Omega,& t>0,\\
u_{t}&+&\nabla P&=\Delta u+n\nabla\phi,\ &x\in\Omega,& t>0,\\
&&\nabla\cdot u&=0,\ &x\in\Omega,& t>0,
\end{array}\right.
\end{align*}
under homogeneous Neumann boundary conditions for $n$ and $c$ and homogeneous Dirichlet boundary conditions for $u$, where $\Omega\subset\mathbb{R}^2$ is a bounded domain with smooth boundary and $\phi\in C^{2}\left(\bar{\Omega}\right)$. From recent results it is known that for suitable regular initial data, the corresponding initial-boundary value problem possesses a global generalized solution. We will show that for small initial mass $\int_{\Omega}\!n_0$ these generalized solutions will eventually become classical solutions of the system and obey certain asymptotic properties.

Moreover, from the analysis of certain energy-type inequalities arising during the investigation of the eventual regularity, we will also derive a result on global existence of classical solutions under assumption of certain smallness conditions on the size of $n_0$ in $L^1\!\left(\Omega\right)$ and in $L\log L\!\left(\Omega\right)$, $u_0$ in $L^4\!\left(\Omega\right)$, and of $\nabla c_0$ in $L^2\!\left(\Omega\right)$.
}

{\noindent\textbf{Keywords:} chemotaxis, Stokes, chemotaxis-fluid interaction, global existence, generalized solution, eventual regularity, stabilization}

{\noindent\textbf{MSC (2010):} 35B65, 35B40 (primary), 35K35, 35Q92, 92C17}
\end{abstract}

\newpage
\section{Introduction}\label{sec1:intro}
Even among the smallest and most primitive organisms there are cases of complex and macroscopical collective behavior, for instance bacteria of species \emph{E. coli} were confirmed to form migrating bands when subjected to a test environment featuring gradients of nutrient concentration (\cite{Adler708}). Following these experimental findings, chemotaxis systems with singular sensitivity of the form
\begin{align}\label{ctsing}
\left\{
\begin{array}{r@{\,}l@{\quad}l}
n_{t}&=\Delta n-\nabla\!\cdot(\frac{n}{c}\nabla c),\\
c_{t}&=\Delta c-nc,\\
\end{array}\right.
\end{align}
were among the first phenomenological models proposed by Keller and Segel (\cite{KS71travbands}) to study these processes of chemotactic migration. Herein, $n$ denotes the density of the bacteria which orient their movement towards increasing concentration $c$ of a chemical substance which serves as their food source and is thereby consumed in the process.
Singular chemotactic sensitivities of the type featured in \eqref{ctsing} express the system assumption that the signal is perceived as described by the Weber-Fechner law (\cite{HP09},\cite{ROSEN1978}).
An outstanding facet of this system, as already illustrated in \cite{KS71travbands}, is the occurence of wave-like solution behavior without any type of cell kinetics, which is known to be vital for such effects in standard reaction-diffusion equations. For studies on existence and stability properties of traveling wave solutions of \eqref{ctsing} see \cite{Wang13surv,LiLiWang14,NagIke91} and references therein.

The results on global existence to systems of the form \eqref{ctsing} are very sparse, with widely arbitrary initial data only being treated for the one-dimensional case (\cite{TWWDCDSB13},\cite{li2015initial}). In higher dimensions the results were constrained to the Cauchy problem for \eqref{ctsing} in $\R^n$ with $n\in\{2,3\}$, where smallness conditions on the initial data had to be imposed to show the existence of globally defined classical solutions (\cite{Wang20162225}). Only recently (\cite{Win16CS1}), so called global generalized solutions to \eqref{ctsing} were constructed in the two-dimensional case. The solutions are obtained through the study of a suitably chosen regularization guaranteeing that the regularized chemical concentration is strictly bounded away from zero for all times. These generalized solutions comply with the classical solution concept in the sense that generalized solutions which are sufficiently smooth also solve the system in the classical sense. In a sequel to the previously mentioned work the author furthermore proved that if the initial mass is small these generalized solutions eventually become classical solutions after some (possibly large) waiting time and that the solutions satisfy certain kind of asymptotic properties (\cite{Win16CS2}).
\\[0.1cm]
\noindent{\textbf{Eventual regularity and fluid interaction.} Our interest slightly differing from the system proposed by Keller and Segel, where the model assumes no interaction between bacteria and surroundings, we will consider the case that the bacteria may be affected by their liquid environment. Here, we do not only assume that this interaction occurs by means of transport, but also in form of a feedback between the cells and the fluid velocity stemming from a buoyancy effect assumed in the model development featured in \cite{tuval2005bacterial}. The experimental evidence reported in the latter reference suggests that the chemotactic motion inside the liquid can be substantially influenced by the feedback between cells and fluid, with turbulence emerging spontaneously in population of aerobic bacteria suspended in sessile drops of water. 
As a prototypical model for the description of this phenomenom a system of the form
\begin{align}\label{CNS}
\left\{
\begin{array}{r@{\,}c@{\,}c@{\,}l}
n_{t}&+&u\cdot\!\nabla n\ &=\Delta n-\nabla\!\cdot(n\nabla c),\\
c_{t}&+&u\cdot\!\nabla c\ &=\Delta c-nc,\\
u_{t}&+&\kappa(u\cdot\nabla)u\ &=\Delta u+n\nabla\phi-\nabla P\\
&&\dive u\ &=0,
\end{array}\right.
\end{align}
was proposed in \cite{tuval2005bacterial} and has been the groundwork for many articles concerning the mathematical analysis of chemotaxis-fluid interaction since the first analytical results asserting local existence of weak solutions (\cite{lorz10}). Obtaining results concerning the global existence of solutions is far from trivial, even when $u\equiv0$ the global existence of solutions is only known under a smallness condition on the initial data (\cite{Tao-consumption_JMAA11}), or when $N=2$ (e.g. \cite{win_fluid_final}). These outcomes are similar in the case of $u\not\equiv0$. In the two-dimensional setting global classical solutions stemming from reasonably smooth initial data have also been shown to exist in \cite{win_fluid_final}, whereas many results treating variants of \eqref{CNS} in three-dimensional frameworks are again restricted to weak solutions emanating from small initial data (e.g. \cite{kozono15},\cite{caolan16_smalldatasol3dnavstokes}). Nevertheless, even in theses cases, where global regularity is hard to prove, some results concerning eventual regularity of solutions have been shown. In particular, for the fluid free case eventual smoothness of solutions was shown in \cite{TaoWin-evsmooth_JDE12} for $N=3$ and a result including fluid is contained in \cite{win15_chemonavstokesfinal}, where certain weak eventual energy solutions are considered.

Similar smoothing effects can also be observed in a setting where $N=3$ and logistic growth terms of the form $+\rho n-\mu n^2$ $(\rho\geq0,\mu>0)$ are included in the first equation. In this framework it is still unclear whether global classical solutions exist for small $\mu>0$ and reasonably arbitrary initial data, but weak solutions which eventually become smooth are known to exist for any $\mu>0$ and possibly large initial data, as indicated by the studies in e.g. \cite{Lan16_M3AS}.

\noindent{\textbf{Chemotaxis-fluid system with singular sensitivity.} 
In light of the regularizing effects observed in the chemotaxis and chemotaxis-fluid problems mentioned above it seems reasonable to assume that also in the case of singular sensitivity the smoothing effect of the second equation will eventually result in classical solutions even if fluid interaction with the bacteria is present. As the construction of weak solution used in \cite{Wang2016} does not work for the full Navier-Stokes subsystem (as included in \eqref{CNS}) we instead work with the simpler Stokes realization of the fluid, which was also employed in \cite{Wang2016}, instead. In fact we will study systems of the form
\begin{align}\label{CN}
\left\{
\begin{array}{r@{\,}c@{\,}c@{\,}l@{\quad}l@{\quad}l@{\,}c}
n_{t}&+&u\cdot\!\nabla n\ &=\Delta n-\nabla\!\cdot(\frac{n}{c}\nabla c),\ &x\in\Omega,& t>0,\\
c_{t}&+&u\cdot\!\nabla c\ &=\Delta c-nc,\ &x\in\Omega,& t>0,\\
u_{t}&+&\nabla P\ &=\Delta u+n\nabla\phi,\ &x\in\Omega,& t>0,\\
&&\dive u\ &=0,\ &x\in\Omega,& t>0,
\end{array}\right.
\end{align}
with boundary conditions
\begin{align}\label{BC}
\frac{\partial n}{\partial\nu}=\frac{\partial c}{\partial \nu}=0,\quad\text{and}\quad u=0\quad\text{for }x\in\romega\text{ and }t>0,
\end{align}
and initial conditions
\begin{align}\label{IC}
n(x,0)=n_0(x),\quad c(x,0)=c_0(x),\quad u(x,0)=u_0(x),\quad x\in\Omega.
\end{align}
$\Omega\subset\R^2$ denotes a bounded domain with smooth boundary and the gravitational potential $\phi$ is assumed to satisfy
\begin{align}\label{phireg}
\phi\in\CSp{2}{\bomega}\quad\text{with}\quad K_1:=\|\phi\|_{\W[1,\infty]}.
\end{align}
For the initial distributions we will prescribe the regularity assumptions
\begin{align}\label{IR}
\left\{\begin{array}{r@{\,}l@{\quad}l}
n_0&\in\CSp{0}{\bomega}&\text{with } n_0\geq0\text{ in }\Omega\text{ and }n_0\not\equiv0,\\
c_0&\in\W[1,\infty]&\text{with }c_0>0\text{ in }\bomega,\\ 
u_0&\in\DAr{\alpha}{r}&\text{for all }r\in(1,\infty)\text{ and some }\alpha\in(\frac12,1),
\end{array}\right.
\end{align}
with $A_r$ denoting the Stokes operator $A_r:=-\HPr\Delta$ in $\LSp{r}{\Omega;\R^2}$ with domain $D\left(A_r\right)=\WSp{2,r}{\Omega;\R^2}\cap W^{1,r}_0\!\left(\Omega;\R^2\right)\cap L^r_\sigma\!\left(\Omega\right)$, where $L^r_\sigma\!\left(\Omega\right)=\{\varphi\in\LSp{r}{\Omega;\R^2}\,\vert\,\dive\varphi=0\}$ stands for the solenodial subspace of $\LSp{r}{\Omega,\R^2}$ obtained by the Helmholtz projection $\HPr$. 

In this setting, building on the work \cite{Win16CS1}, it was shown in \cite{Wang2016} that for any $(n_0,c_0,u_0)$ satisfying \eqref{IR} the system \eqref{CN} possesses at least on global generalized solution (in the sense of Definition \ref{def:gen_sol} below). These solutions are constructed by a similar limiting procedure as in the fluid free setting, making sure that for each of the approximate solutions the quantity $c$ remains strictly positive throughout $\Omega$ for all times. In a simplified version the result on global existence of generalized solutions and basic decay properties of $c$ obtained in \cite{Wang2016} can be summarized as follows.

\begin{thmx}\label{thm:globsol}\ \\
Let $\Omega\subset\R^2$ be a bounded domain with smooth boundary. Then for all $(n_0,c_0,u_0)$ satisfying \eqref{IR}, the problem \eqref{CN}--\,\eqref{IC} possesses at least one global generalized solution $(n,c,u)$ in the sense of Definition \ref{def:gen_sol} below. For each $p\in[1,\infty)$ the solution has the properties that $n(\cdot,t)\in\Lo[p]$ and $\frac{\nabla c}{c}\in\Lo[2]$ for a.e. $t>0$. Moreover, $c$ is continuous on $[0,\infty)$ as $\Lo[\infty]$--valued function with respect to the weak--$\star$ topology on $\Lo[\infty]$, and satisfies
\begin{align*}
c(\cdot,t)\wsto0\quad\text{in }\Lo[\infty]\qquad\text{and}\qquad c(\cdot,t)\to0\quad\text{in }\Lo[p]\qquad\text{as }t\to\infty.
\end{align*}
\end{thmx}

{\noindent \textbf{Main results.}} The existence of global generalized solutions as provided by Theorem \ref{thm:globsol} at hand, it is the purpose of the present work to study the question how far the eventual regularity and stabilization results for small data, as obtained in \cite{Win16CS2} for \eqref{ctsing}, may be affected by the interaction of the bacteria with their liquid surroundings.

\begin{theorem}\label{thm:evsmooth}
Let $\Omega\subset\R^2$ be a bounded domain with smooth boundary. Then there exists some $\ms\!>0$ such that for any $(n_0,c_0,u_0)$ satisfying \eqref{IR} as well as
\begin{align}\label{eq:critmass}
\intomega n_0\leq\ms,
\end{align}
the global generalized solution of \eqref{CN}--\,\eqref{IC} from Theorem \ref{thm:globsol} has the property that there exists $T>0$ such that
\begin{align}\label{eq:evreg}
n\in\CSp{2,1}{\bomega\times[T,\infty)},\quad c\in\CSp{2,1}{\bomega\times[T,\infty)}\quad\text{and}\quad u\in\CSp{2,1}{\bomega\times[T,\infty);\R^2},
\end{align}
that
\begin{align}\label{eq:large-time-positivity-c}
c(x,t)>0\quad\text{for all }x\in\bomega\text{ and any }t\geq T,
\end{align}
and such that $(n,c,u)$ solve \eqref{CN}--\,\eqref{IC} classically in $\Omega\times(T,\infty)$. Furthermore, this solution satisfies
\begin{align}\label{eq:conv-n}
n(\cdot,t)&\to \frac{1}{|\Omega|}\intomega n_0\quad\text{in }\Lo[\infty],
\qquad c(\cdot,t)\to 0\quad\ \text{in }\Lo[\infty],\qquad
u(\cdot,t)\to 0\quad\ \text{in }\Lo[\infty],
\end{align}
and
\begin{align}\label{eq:conv-nabc}
\frac{\nabla c(\cdot,t)}{c(\cdot,t)}\to0\quad\text{in }\Lo[\infty]
\end{align}
as $t\to\infty$.
\end{theorem}

Our analysis will also in straightforward manner allow us to formulate a result for global classical solutions to \eqref{CN}--\,\eqref{IC} if certain smallness conditions are fulfilled by the initial distributions. Furthermore, these global classical solutions inherit the same asymptotic properties stated in Theorem \ref{thm:evsmooth}. In order to completely formulate this outcome, we note that in two-dimensional domains by the \GNI\ and elliptic regularity theory one can find $K_2>0$ and $K_3>0$ such that
\begin{align}\label{global_gnb_constants1} 
\|\varphi\|_{\Lo[3]}^3&\leq K_2\|\varphi\|_{\W[1,2]}^2\|\varphi\|_{\Lo[1]}\quad\text{for all }\varphi\in \W[1,2]
\end{align}
and
\begin{align}\label{global_gnb_constants3}
\|\nabla\varphi\|_{\Lo[4]}&\leq K_3\|\Delta\varphi\|_{\Lo[2]}^\nfrac{1}{2}\|\nabla\varphi\|_{\Lo[2]}^\nfrac{1}{2}\quad\text{for all }\varphi\in\W[2,2]\text{ with }\frac{\partial\varphi}{\partial\nu}=0\text{ on }\romega.
\end{align}
We obtain the following.
\begin{theorem}\label{thm:smalldataglobalclass}
Let $\Omega\subset\R^2$ be a bounded domain with smooth boundary. Then there exists $\mss>0$ such that
for any $(n_0,c_0,u_0)$ satisfying \eqref{IR},
\begin{align}\label{eq:init-small}
\intomega n_0\leq \mss,\quad\text{and}\quad \intomega|u_0|^4\leq\mss
\end{align}
as well as
\begin{align}\label{eq:init-energ}
\intomega n_0\ln\frac{n_0}{\mu}+\frac{1}{2}\intomega\frac{|\nabla c_0|^2}{c_0^2}<\min\left\{\frac{1}{4K_3},\frac{1}{8K_2}\right\}-\frac{\mu|\Omega|}{e}
\end{align}
for some $\mu>0$ and $K_2$, $K_3$ given by \eqref{global_gnb_constants1} and \eqref{global_gnb_constants3}, repsectively, there exists a triple $(n,c,u)$ of functions, for each $\vartheta>2$ uniquely determined by the inclusions
\begin{align*}
\begin{cases}
n\in\CSp{0}{\bomega\times[0,\infty)}\cap\CSp{2,1}{\bomega\times(0,\infty)},\\
c\in\CSp{0}{\bomega\times[0,\infty)}\cap\CSp{2,1}{\bomega\times(0,\infty)}\cap\LSploc{\infty}{[0,\infty);\W[1,\vartheta]},\\
u\in\CSp{0}{\bomega\times[0,\infty);\R^2}\cap\CSp{2,1}{\bomega\times(0,\infty);\R^2},\end{cases}
\end{align*}
such that $n>0$ in $\bomega\times(0,\infty)$ and $c>0$ in $\bomega\times[0,\infty)$, and such that $(n,c,u)$ together with some $P\in\CSp{1,0}{\bomega\times[0,\infty)}$ solve \eqref{CN}--\,\eqref{IC} in the classical sense in $\Omega\times(0,\infty)$. Furthermore, this solution has the convergence properties stated in Theorem \ref{thm:evsmooth}.
\end{theorem}

In contrast to the known result for the system without fluid, obtained by taking $u\equiv0$ in \eqref{CN} where requiring only $\intomega n_0\ln\frac{n_0}{\mu}+\frac{1}{2}\intomega|\frac{|\nabla c_0|^2}{c_0^2}$ to be small was sufficient to obtain global classical solutions, in this case we require  additional smallness conditions in the form of sufficiently small bounds for $n_0$ in $\Lo[1]$ and $u_0$ in $\Lo[4]$.
\\[0.1cm]
{\noindent \textbf{Notation.}} Throughout the article, in addition to the previously mentioned assumptions in \eqref{phireg} and \eqref{IR} for $\Omega$, $\phi$, the initial data, the Stokes operator and its semigroup, we will make use of the following notations. $\lambda_1>0$ will always denote the first positive eigenvalue of the Stokes operator in $\Omega$ with respect to homogeneous Dirichlet boundary data. Since $A^\alpha_r\varphi,e^{-tA_r}\varphi$ and $\HPr\psi$ are independent of $r\in(1,\infty)$ for $\varphi\in C_0^\infty\!\left(\Omega\right)\cap L_\sigma^r\!\left(\Omega\right)$ and $\psi\in C_0^\infty\!\left(\Omega\right)$, we will drop the subscript whenever there is no danger of confusion. Similar to denoting by $L_\sigma^r\!\left(\Omega\right)$ all divergence free functions of $\Lo[p]$, the space of divergence free, smooth test functions with compact support in $\Omega\times(0,\infty)$ will be denoted by $C^\infty_{0,\sigma}\!\left(\Omega\times(0,\infty)\right)$. Additionally, when talking about classical solutions to some of the featured systems in $\Omega\times(t_0,\infty)$ for some $t_0\geq0$, we will often shorten the notation to $(n,c,u)\in\CSp{0}{\Omega\times[t_0,\infty)}$, when we are actually considering $(n,c,u,P)\in\CSp{0}{\Omega\times[t_0,\infty)}\times\CSp{0}{\Omega\times[t_0,\infty)}\times\CSp{0}{\Omega\times[t_0,\infty);\R^2}\times\CSp{1,0}{\bomega\times[t_0,\infty)}$. The notation $(n,c,u)\in\CSp{2,1}{\Omega\times(t_0,\infty)}$ will be used in a similar fashion.

\setcounter{equation}{0} 
\section{Basic properties of a family of generalized problems}\label{sec2:approxsol}
The construction of the generalized solution mentioned above is based on a limit procedure of solutions to regularized problems and a transformation thereof. Since the original problem \eqref{CN} and the family of approximate problems in question are quite similar, we will first consider the even more general family of problems
\begin{align}\label{CSc}
\arraycolsep=1.4pt\def\arraystretch{1.25}
\left\{
\begin{array}{r@{\,}c@{\,}c@{\,}l@{\quad}l@{\quad}l@{\,}c}
n_{t}&+&u\cdot\!\nabla n\ &=\Delta n-\nabla\!\cdot\Big(\frac{nf'(n)}{c}\nabla c\Big),\ &x\in\Omega,& t>0\\
c_{t}&+&u\cdot\!\nabla c\ &=\Delta c-f(n)c,\ &x\in\Omega,& t>0,\\
u_{t}&+&\nabla P &=\Delta u+n\nabla\phi,\ &x\in\Omega,& t>0,\\
&&\dive u\ &=0,\ &x\in\Omega,& t>0,
\end{array}\right.
\end{align}
where we only require that the functions $f\in\CSp{3}{[0,\infty)}$ satisfy
\begin{align}\label{eq:fnoepsprop}
f(0)=0\quad\text{ and }\quad0\leq f'\leq 1\text{ on }[0,\infty).
\end{align}
Upon proper choice of a subfamily of these functions (c.f. \eqref{eq:feps_def} below) the system will be regularized in a way that ensures that $c$ is bounded away from zero, from which one can easily obtain global and bounded solutions to the corresponding approximate problems. These global and bounded solutions are one of the main ingredients of the limit process involved in the construction of the generalized solution (\cite{Win16CS1},\cite{Wang2016}).

The problems \eqref{CSc} will be regarded under the boundary conditions
\begin{align}\label{BCc}
\frac{\partial n}{\partial\nu}=\frac{\partial c}{\partial \nu}=0,\quad\text{and}\quad u=0\quad\text{for }x\in\romega\text{ and }t\in(0,\Tm),
\end{align}
and the initial conditions
\begin{align}\label{ICc}
n(x,0)=n_0(x),\quad c(x,0)=c_0(x),\quad u(x,0)=u_0(x),\quad x\in\Omega.
\end{align}
For any $f\in\CSp{3}{[0,\infty)}$ satisfying the conditions above, local existence of classical solutions can be obtained by well-established fixed point methods. Since the necessary adaptions are quite straightforward, we will refer to local existence proofs in closely related situations for details.

\begin{lemma}\label{lem:loc_ex_approx_c}
Let $\Omega\subset\R^2$ be a bounded domain with smooth boundary, $\vartheta>2$ and $f\in\CSp{3}{[0,\infty)}$ satisfies \eqref{eq:fnoepsprop}. Then for all $(n_0,c_0,u_0)$ satisfying \eqref{IR} there exist $\Tm\in(0,\infty]$ and uniquely determined functions
\begin{align*}
n&\in\CSp{0}{\bomega\times[0,\Tm)}\cap\CSp{2,1}{\bomega\times(0,\Tm)},\\
c&\in\CSp{0}{\bomega\times[0,\Tm)}\cap\CSp{2,1}{\bomega\times(0,\Tm)}\cap\CSp{0}{[0,\Tm);\W[1,\vartheta]},\\
u&\in\CSp{0}{\bomega\times[0,\Tm);\R^2}\cap\CSp{2,1}{\bomega\times(0,\Tm);\R^2},
\end{align*}
which together with some $P\in\CSp{1,0}{\bomega\times[0,\Tm)}$ solve \eqref{CSc}--\,\eqref{ICc} in the classical sense and satisfy $n>0$ and $c>0$ in $\bomega\times(0,\Tm)$ as well as 
\begin{align}\label{eq:loc_ex_alt}
\quad\Tm=\infty,\,&\mbox{ or }\,\liminf_{t\nearrow\Tm}\inf_{x\in\Omega}c(x,t)=0,\\&\mbox{ or }\,\limsup_{t\nearrow\Tm}\left(\|n(\cdot,t)\|_{\Lo[\infty]}+\|c(\cdot,t)\|_{\W[1,\vartheta]}+\|A^{\alpha}u(\cdot,t)\|_{\Lo[2]}\right)=\infty.\nonumber
\end{align}
Furthermore, the solution has the properties that
\begin{align}\label{eq:mass_cons_n}
\intomega n(x,t)\intd x=\intomega n_0(x)\intd x\quad\text{for all }t\in(0,\Tm)
\end{align}
and
\begin{align}\label{eq:c_bound_globex}
&c(x,t)\leq\|c_0\|_{\Lo[\infty]}\quad\text{for all }(x,t)\in\bomega\times[0,\Tm).
\end{align}
\end{lemma}

\begin{bew}
Local existence, uniqueness and the blow-up criterion \eqref{eq:loc_ex_alt} can be obtained by straightforward adaption of well known arguments as detailed in \cite{HoWin05_bvblowchemo,FIY14} and \cite{win_fluid_final} for related situations. Simple integration of the first equation in \eqref{CSc} proves \eqref{eq:mass_cons_n}, whereas by the nonnegativity of $f$ an application of the parabolic comparison principle to the second equation in \eqref{CSc}, with $\bar{c}\equiv\|c_0\|_{\Lo[\infty]}$ taken as supersolution, immediately entails \eqref{eq:c_bound_globex}.
\end{bew}

\subsection{Regularity of the Stokes subsystem}\label{sec2:regu}
It is well known that the Stokes subsystem $\frac{\intd}{\intd t}u+ A u=\HP(n\nabla\phi)$ in \eqref{CSz} has the property that the regularity of the spatial derivative $\nabla u$ is solely reliant on the regularity of $n$ (since $\nabla \phi$ is bounded). In fact for Stokes systems of the form
\begin{align}\label{stokessys}
\arraycolsep=1.4pt\def\arraystretch{1.25}
\left\{
\begin{array}{r@{\,}l@{\quad}l@{\quad}l@{\,}c}
u_{t}&=\Delta u-\nabla P+n\nabla\phi,\ &x\in\Omega,& t_0>0,\\
\dive u\ &=0,\ &x\in\Omega,& t_0>0,\\
u\ &=0,\ &x\in\romega,& t_0>0,
\end{array}\right.
\end{align}
we can obtain the following two results. The first is a refinement of a basic boundedness result e.g. featured in \cite[Lemma 2.4]{Wang20157578}.

\begin{lemma}\label{lem:l4_bound_u}
Let $\phi\in\CSp{2}{\bomega}$. There exist constants $\lambda_1>0$ and $K_{u}>0$ such that whenever 
$u\in\CSp{0}{\bomega\times[t_0,T_0);\R^2}\cap\CSp{2,1}{\bomega\times(t_0,T_0);\R^2}$ is a classical solution of \eqref{stokessys} in $\Omega\times(t_0,T_0)$ for some $0\leq t_0<T_0\leq\infty$ and $n\in\CSp{0}{\bomega\times[t_0,T_0)}$ satisfying
\begin{align*}
\intomega |n(\cdot,t)|\leq L\quad\text{for all }t\in(t_0,T_0),
\end{align*}
with some $L>0$, then
\begin{align*}
\|u(\cdot,t)\|_{\Lo[4]}&\leq K_{u}e^{-\lambda_1 (t-t_0)}\|u(\cdot,t_0)\|_{\Lo[4]}+K_{u}L\quad\text{for all }t\in(t_0,T_0).
\end{align*}
\end{lemma}

\begin{bew}
By the variation-of-constants representation for $u$ we have
\begin{align*}
u(\cdot,t)=e^{-(t-t_0)A}u(\cdot,t_0)+\int_{t_0}^te^{-(t-s)A}\HP(n(\cdot,s)\nabla\phi)\intd s\quad\text{for all }t\in(t_0,T_0).
\end{align*}
Fixing any $\gamma\in(\frac{3}{4},1)$ we see that
\begin{align*}
\|u(\cdot,t)\|_{\Lo[4]}\leq \|e^{-(t-t_0)A}u(\cdot,t_0)\|_{\Lo[4]}+\int_{t_0}^t\|A^\gamma e^{-(t-s)A}A^{-\gamma}\HP(n(\cdot,s)\nabla\phi)\|_{\Lo[4]}\intd s
\end{align*}
holds for all $t\in(t_0,T_0)$. Now, in view of the well known regularity estimates for the Stokes semigroup (e.g. \cite[Lemma 3.1]{win_ct_fluid_3d}) we find constants $\lambda_1>0$ and $C_1>0$  such that
\begin{align*}
\|e^{-(t-t_0)A}u(\cdot,t_0)\|_{\Lo[4]}\leq C_1 e^{-\lambda_1 (t-t_0)}\|u(\cdot,t_0)\|_{\Lo[4]}\quad\text{for all }t>t_0,
\end{align*}
and, since for $1\leq p<q<\infty$ and $\gamma\in(0,1)$ satisfying $\gamma>\frac{1}{p}-\frac{1}{q}$ it holds that $\|A^{-\gamma}\HP\varphi\|_{\Lo[q]}\leq C\|\varphi\|_{\Lo[p]}$ for all $\varphi\in C_0^\infty\!\left(\Omega\right)$ (\cite[Lemma 2.3]{Wang20157578}), there exists $C_2>0$ such that
\begin{align*}
\|A^\gamma e^{-(t-s)A}A^{-\gamma}\HP(n(\cdot,s)\nabla\phi)\|_{\Lo[4]}
\leq C_2(t-s)^{-\gamma}e^{-\lambda_1 (t-s)}\|n(\cdot,s)\nabla\phi\|_{\Lo[1]}\quad\text{for all }s\in(t_0,t),
\end{align*}
by choice of $\gamma\in(\frac{3}{4},1)$. Hence, relying on \eqref{phireg} and our assumption for $\intomega |n(\cdot,t)|$, we may estimate
\begin{align*}
\|u(\cdot,t)\|_{\Lo[4]}\leq C_1e^{-\lambda_1 (t-t_0)}\|u(\cdot,t_0)\|_{\Lo[4]}+C_2K_1L\int_0^\infty(t-s)^{-\gamma}e^{-\lambda_1 (t-s)}\intd s\quad\text{for all }t>0,
\end{align*}
which due to $\gamma<1$ concludes the proof upon obvious choice for $K_u$.
\end{bew}

The second lemma regarding the Stokes subsystem concerns norms of the spatial gradient of $u$. These results are well-known. (see e.g. \cite[Lemma 2.5]{Wang20157578} and \cite[Corollary 3.4]{win_ct_fluid_3d} for details.)

\begin{lemma}\label{lem:stokes_est}
Assume $\alpha\in(\frac{1}{2},1)$, $t_0\geq0$ and $\phi\in\CSp{2}{\bomega}$ and let $p\in[1,\infty)$ and $r\in[1,\infty]$ be such that
\begin{align*}
\begin{cases}
r<\frac{2p}{2-p}\quad&\text{if }p\leq2,\\
r\leq\infty&\text{if }p>2.
\end{cases}
\end{align*}
Then for any $u(\cdot,t_0)\in\DAr{\alpha}{r}$ there exists a constant $C=C(u(\cdot,t_0),\phi,p,r,L)>0$ such that whenever 
$u\in\CSp{0}{\bomega\times[t_0,T_0)}\cap\CSp{2,1}{\bomega\times(t_0,T_0)}$ is a classical solution of \eqref{stokessys}
 in $\Omega\times(t_0,T_)$ for some $0\leq t_0<T_0\leq\infty$ and $n\in\CSp{0}{\bomega\times[t_0,T)}$ satisfying
\begin{align*}
\|n(\cdot,t)\|_{\Lo[p]}\leq L\quad\text{for all }t\in(t_0,T),
\end{align*}
with some $L>0$, then 
\begin{align*}
\|\nabla u(\cdot,t)\|_{\Lo[r]}\leq Ce^{-\lambda_1 (t-t_0)}+CL\quad\text{for all }t\in(t_0,T).
\end{align*}
\end{lemma}

In particular, in view of the mass conservation property of $n$ and the Sobolev embedding theorem, we can easily obtain bounds independent of $f$ for the quantity $\|u\|_{\Lo[p]}$ with $p<\infty$ from the previous Lemma. For these potentially better bounds than the one provided by Lemma \ref{lem:l4_bound_u} however, we do not know the exact relation to $u_0$.

\subsection{Logarithmic rescaling and basic a priori information on \texorpdfstring{$z$}{z}}\label{sec2:trafo}
Now, a quite standard change in variables transformation obtained by taking $n$, $c$ and $u$ from Lemma \ref{lem:loc_ex_approx_c} and setting
\begin{align*}
z:=-\ln\left(\frac{c}{\|c_0\|_{\Lo[\infty]}}\right)\quad\text{ and }\quad z_0:=-\ln\left(\frac{c_0}{\|c_0\|_{\Lo[\infty]}}\right),
\end{align*}
will lead to the transformed systems
\begin{align}\label{CSznoeps}
\left\{
\begin{array}{r@{\,}c@{\,}c@{\,}l@{\quad}l@{\quad}l@{\,}c}
n_{t}&+&u\cdot\!\nabla n\ &=\Delta n+\nabla\!\cdot(n f'(n)\nabla z),\ &x\in\Omega,& t>0,\\
z_{t}&-&u\cdot\!\nabla z\ &=\Delta z-|\nabla z|^2+f(n),\ &x\in\Omega,& t>0,\\
u_{t}&+&\nabla P &=\Delta u+n\nabla\phi,\ &x\in\Omega,& t>0,\\
&&\dive u\ &=0,\ &x\in\Omega,& t>0,
\end{array}\right.
\end{align}
which build the basis for our analysis of the energy-type inequalities featured in Section \ref{sec4:energyfunctional}. This transformation has been thoroughly used in previous literature (see e.g. \cite{Wang20162225},\cite{Win16CS1},\cite{Win16CS2}) to analyze systems in similar settings. We will consider \eqref{CSznoeps} along with the boundary conditions
\begin{align}\label{CSznoepsBC}
\frac{\partial n}{\partial\nu}=\frac{\partial z}{\partial \nu}=0,\quad\text{and}\quad u=0\quad\text{for }x\in\romega\text{ and }t>0,
\end{align}
and initial conditions
\begin{align*}
n(x,0)=n_0(x),\quad z(x,0)=z_0(x):=-\ln\left(\frac{c_0(x)}{\|c_0\|_{\Lo[\infty]}}\right),\quad u=u_0(x),\quad x\in\Omega.
\end{align*}

\begin{remark}\label{rem:masscons-sys-csz}
Let $f\in\CSp{3}{[0,\infty)}$ satisfy \eqref{eq:fnoepsprop}. Assume that $(n,z,u)\in\CSp{2,1}{\bomega\times(T_1,T_2)}$ is a classical solution of the boundary value problem \eqref{CSznoeps},\eqref{CSznoepsBC} in $\Omega\times(T_1,T_2)$ with some $T_1\geq0$ and $T_2\in(T_1,\infty]$. Then the solution satisfies the mass conservation property
\begin{align*}
\frac{\intd}{\intd t}\intomega n(\cdot,t)=0\quad\text{for all }t\in(T_1,T_2).
\end{align*}
\end{remark}
This reformulation of our previous generalized systems at hand, we immediately obtain the following basic information -- not depending on $f$ -- about the transformed chemical concentration $z$.

\begin{lemma}\label{lem:bound_zeps}
Let $m_0>0$. Suppose that for $f\in C^3\!\left([0,\infty)\right)$ satisfying \eqref{eq:fnoepsprop} and $t_0\geq0$ the triple $(n,z,u)\in\CSp{2,1}{\bomega\times(t_0,\infty)}$ is a classical solution of \eqref{CSznoeps}--\eqref{CSznoepsBC} in $\Omega\times(t_0,\infty)$ with the properties that $n\geq0$ in $\Omega\times(t_0,\infty)$ and $\intomega n(\cdot,t_0)\leq m_0$. Then
\begin{align}\label{eq:bound_zeps}
\intomega z(\cdot,t)+\int_{t_0}^t\!\intomega |\nabla z|^2\leq \intomega z(\cdot,t_0)+(t-t_0)m_0\quad\text{for all }t>t_0.
\end{align}
\end{lemma}

\begin{bew}
Integrating the second equation of \eqref{CSznoeps} with respect to space shows that
\begin{align*}
\frac{\intd}{\intd t}\intomega z=\intomega\Delta z-\intomega|\nabla z|^2+\intomega f(n)+\intomega u\cdot\nabla z
\end{align*}
holds for all $t\in(t_0,\infty)$. Making use of $\dive u=0$, the Neumann boundary conditions for $z$, $n\geq0$ and the fact that $f(s)\leq s$ for all $s\geq0$ we obtain, upon integration by parts, that
\begin{align*}
\frac{\intd}{\intd t}\intomega z+\intomega|\nabla z|^2\leq\intomega n
\end{align*}
is valid on $t\in(0,\infty)$. Due to the mass conservation we have $\intomega n(\cdot,t)\leq m_0$ for all $t>t_0$ and therefore integrating this inequality immediately establishes \eqref{eq:bound_zeps}.
\end{bew}

\setcounter{equation}{0} 
\section{Generalized solution concept and approximate solutions}\label{sec3:uandz}
Before going into more detail for our eventual smoothness result, let us briefly review the solution concept of generalized solutions and the exact form of the approximate problems. These were already used in \cite{win15_chemorot,Win16CS1} for the closely related settings without fluid and in \cite{Wang2016} for the system with Stokes fluid.

A global generalized solution is defined as follows (see also \cite[Definition 2.1--2.3]{win15_chemorot},\cite[Definition 2.1]{Wang2016}).

\begin{definition}\label{def:gen_sol}
Assume that $(n_0,c_0,u_0)$ satisfy \eqref{IR}. Suppose that a triple $(n,c,u)$ of functions
\begin{align}\label{eq:sol_def_reg1}
\arraycolsep=1.4pt\def\arraystretch{1.25}
\left\{\begin{array}{l@{\, }c@{\, }l}
n&\in&\LSploc{1}{\bomega\times[0,\infty)},\\
c&\in&\LSploc{\infty}{\bomega\times[0,\infty)}\cap\LSploc{2}{[0,\infty);\W[1,2]},\\
u&\in& L^{1}_{loc}\big([0,\infty); W^{1,1}_{0}\left(\Omega\right);\R^2\big),
\end{array}\right.
\end{align}
satisfies
\begin{align}\label{eq:sol_def_pos}
n\geq0,\quad\text{and}\quad c>0,\quad\text{and}\quad\dive u=0\quad\text{a.e. in }\Omega\times(0,\infty),
 \end{align}
as well as
\begin{align}\label{eq:sol_def_reg2}
\nabla\ln(n+1)\in\LSploc{2}{\bomega\times[0,\infty)}\qquad\text{and}\qquad\nabla\ln c\in\LSploc{2}{\bomega\times[0,\infty)}.
\end{align}
Then $(n,c,u)$ will be called a global generalized solution of \eqref{CN}--\,\eqref{IC} if $n$ satisfies the mass conservation property
\begin{align*}
\intomega n(x,t)\intd x=\intomega n_0(x)\intd x\quad\text{for a.e. }t>0,
\end{align*}
if the inequality
\begin{align}\label{eq:sol_def_nineq}
-&\intinfomega\ln(n+1)\varphi_t-\intomega\ln(n_0+1)\varphi(\cdot,0)\nonumber\\
&\geq\intinfomega|\nabla \ln(n+1)|^2\varphi-\intinfomega\nabla\ln(n+1)\cdot\nabla\varphi+\intinfomega\frac{n}{n+1}\nabla \ln c\cdot\nabla \varphi\\
&\quad-\intinfomega\frac{n}{n+1}\left(\nabla\ln(n+1)\cdot\nabla\ln c\right)\varphi+\intinfomega \ln(n+1)(u\cdot\nabla\varphi)\nonumber
\end{align}
holds for each nonnegative $\varphi\in C^\infty_0\left(\bomega\times[0,\infty)\right)$, if the identity
\begin{align*}
\ &\intinfomega c\psi_t+\intomega c_0\psi(\cdot,0)=\intinfomega\nabla c\cdot\nabla \psi+\intinfomega nc\psi-\intinfomega cu\cdot\nabla\psi
\end{align*}
is valid for any  $\psi\in\LSp{\infty}{\Omega\times(0,\infty)}\cap\LSp{2}{(0,\infty);\W[1,2]}$ compactly supported in $\bomega\times[0,\infty)$ with $\psi_t\in\LSp{2}{\Omega\times(0,\infty)}$, and if furthermore the equality
\begin{align*}
\intinfomega u\cdot\Psi_t+\intomega u_0\cdot\Psi(\cdot,0)=\intinfomega\nabla u\cdot\nabla\Psi-\intinfomega n\nabla\phi\cdot\Psi
\end{align*}
holds for all $\Psi\in C_{0,\sigma}^\infty\left(\Omega\times[0,\infty)\right)$.
\end{definition}

It can easily be verified that the supersolution property in \eqref{eq:sol_def_nineq} combined with the mass conservation \eqref{eq:mass_cons_n} is sufficient to obtain that sufficiently regular global generalized solutions are also global solutions in the classical sense (see \cite[Remark 2.1 ii)]{Win16CS1}), i.e. if $(n,c,u)$ is a global generalized solution in the sense of Definition \ref{def:gen_sol} and satisfies $n\geq0$ and $c>0$ in $\bomega\times[0,\infty)$ as well as $(n,c,u)\in\!\CSp{0}{\bomega\times[0,\infty)}\cap\CSp{2,1}{\bomega\times(0,\infty)}$ then $(n,c,u)$ solves \eqref{CSc} in the classical sense.

Generalized solutions in the sense of Definition \ref{def:gen_sol} are constructed by an approximation procedure relying on regularizations in the form of \eqref{CSznoeps} with suitably chosen $f\equiv f_\eps$ (\cite{Win16CS1,Win16CS2,Wang2016}). For this we first fix a cut-off function $\rho\in\CSp{\infty}{[0,\infty)}$ fulfilling $\rho\equiv1$ in $[0,1]$ and $\rho\equiv0$ in $[2,\infty)$ and define the family of functions $\left\{f_\eps\right\}_{\eps\in(0,1)}\subseteq\CSp{\infty}{[0,\infty)}$ given by
\begin{align}\label{eq:feps_def}
f_\eps(s):=\int_0^s\rho(\eps \sigma)\intd \sigma,\qquad s\geq0.
\end{align}
Every function in this family evidently has the properties
\begin{align}\label{eq:feps_prop1}
f_\eps(0)=0\quad\text{ and }\quad 0\leq f'_\eps\leq 1\quad\text{ on }[0,\infty),
\end{align}
as well as
\begin{align*}
f_\eps(s)=s\quad\text{for all }s\in[0,\textstyle\frac{1}{\eps}]\quad\text{ and }\quad f'_\eps(s)=0\quad\text{for all }s\geq\frac{2}{\eps}.
\end{align*}
Furthermore it holds that
\begin{align*}
f_\eps(s)\nearrow s\quad\text{and}\quad f'_\eps(s)\nearrow 1\quad\text{as }\eps\searrow0
\end{align*}
for each $s\geq0$. According to this choice we can ensure that for the local solutions to \eqref{CSc}\,--\eqref{ICc}  $n_\eps$ is bounded throughout $\Omega\times(0,\Tm)$, and that $c_\eps$ is strictly positive on $\bomega\times(0,\Tm)$, meaning that the most troublesome terms of the extensibility criterion in \eqref{eq:loc_ex_alt} remain bounded, whence the further estimation of remaining less troublesome terms in fact shows that the solution actually is global (\cite{Wang2016}).

Relying on the logarithmic transformation again we obtain for this family of regularizing functions, \eqref{CSznoeps}--\,\eqref{CSznoepsBC} systems of the form

\begin{align}\label{CSz}
\left\{
\begin{array}{r@{\,}c@{\,}c@{\,}l@{\quad}l@{\quad}l@{\,}c}
n_{\eps t}&+&u_{\eps}\cdot\!\nabla n_{\eps}\ &=\Delta n_{\eps}+\nabla\!\cdot(n_{\eps}f'_\eps(n_\eps)\nabla z_{\eps}),\ &x\in\Omega,& t>0,\\
z_{\eps t}&-&u_{\eps}\cdot\!\nabla z_{\eps}\ &=\Delta z_{\eps}-|\nabla z_\eps|^2+f_\eps(n_{\eps}),\ &x\in\Omega,& t>0,\\
u_{\eps t}&+&\nabla P_\eps &=\Delta u_{\eps}+n_{\eps}\nabla\phi,\ &x\in\Omega,& t>0,\\
&&\dive u_{\eps}\ &=0,\ &x\in\Omega,& t>0,
\end{array}\right.
\end{align}
with boundary conditions
\begin{align}\label{BCz}
\frac{\partial n_\eps}{\partial\nu}=\frac{\partial z_\eps}{\partial \nu}=0,\quad\text{and}\quad u_\eps=0\quad\text{for }x\in\romega\text{ and }t>0,
\end{align}
and initial conditions
\begin{align}\label{ICz}
n_\eps(x,0)=n_0(x),\quad z_\eps(x,0)=z_0(x)=-\ln\left(\frac{c_0(x)}{\|c_0\|_{\Lo[\infty]}}\right),\quad u_\eps(x,0)=u_0(x),\quad x\in\Omega.
\end{align}
According to \cite{Wang2016} also these problems posses global classical solutions, with again $n_\eps$ and $z_\eps$ being nonnegative, $n_\eps$ still satisfying the mass conservation property as in Remark \ref{rem:masscons-sys-csz} and $(n_\eps,z_\eps,u_\eps)$ correspond to solutions of systems of the form \eqref{CSc} by means of the substitution $z_\eps=-\ln\left(\frac{c_\eps}{\|c_0\|_{\Lo[\infty]}}\right)$.

The following result summarizes the result on approximation of the generalized solutions established in \cite[Lemma 2.5]{Wang2016}.

\begin{lemma}\label{lem:convergence1}
Let \eqref{IR} hold and denote by $(n,c,u)$ the global generalized solution of \eqref{CN}--\,\eqref{IC} from Theorem \ref{thm:globsol}. Then there exists a sequence $\{\eps_j\}_{j\in\mathbb{N}}\subset(0,1)$ such that $\eps_j\searrow0$ as $j\to\infty$ and such that, for the choice of $f\equiv f_\eps$ in \eqref{CSc}, the corresponding solution $(n_\eps,c_\eps,u_\eps)$ of \eqref{CSc}--\,\eqref{ICc} satisfies
\begin{align*}
n_\eps\to n,\qquad\text{and}\qquad c_\eps\to c,\qquad\text{as well as}\qquad u_\eps\to u\qquad\text{a.e. in }\Omega\times(0,\infty).
\end{align*}
as $\eps=\eps_j\searrow 0$.
\end{lemma}

\setcounter{equation}{0} 
\section{Eventual smoothness of small-data generalized solutions}\label{sec4:smalldatasmooth}
\subsection{Nonincreasing energy for small mass}\label{sec4:energyfunctional}

We will appropriately adjust the functional methods employed in \cite{Win16CS2} to our needs. In fact we will study the behavior of functionals of the form
\begin{align}\label{Fdef}
\F(n,z):=\intomega n\ln\frac{n}{\mu}+\frac{1}{2}\intomega|\nabla z|^2
\end{align}
for $\mu>0$, $0\leq n\in\CSp{0}{\bomega}$ and $z\in\CSp{1}{\bomega}$. We will show that a suitable condition on the size of $\F\big(n(\cdot,t_0),z(\cdot,t_0)\big)$ for some $t_0\geq0$ implies that $F_{\mu}$ is non-increasing from that time onward, along the trajectory of classical solutions to the system \eqref{CSznoeps}. Since we are working with the more generalized version of \eqref{CSz} almost all of the properties of $\F$ also hold in our limit case $f(\xi)\equiv \xi$ obtained by taking $\eps\searrow0$ in \eqref{CSz}. In particular, this will also hold true for the conditional regularity estimates discussed in Section \ref{sec4:regularity}.

We start with some basic relations between $\F$ and the quantities appearing therein.

\begin{lemma}\label{lem:Fprop}
For $\mu>0$ let $\F$ be given by \eqref{Fdef}. Then for all nonnegative $n\in\CSp{0}{\bomega}$ and any $z\in\CSp{1}{\bomega}$ we have
\begin{align}\label{eq:Fineq_n}
\intomega n|\ln n|\leq \F(n,z)+\ln\mu\intomega n+\frac{2|\Omega|}{e},
\end{align}
and
\begin{align}\label{eq:Fineq_z}
\intomega|\nabla z|^2&\leq 2\F(n,z)+\frac{2\mu|\Omega|}{e},
\end{align}
as well as
\begin{align}\label{eq:Fineq_below}
\F(n,z)\geq -\frac{\mu|\Omega|}{e}.
\end{align}
\end{lemma}

\begin{bew}
Making use of the facts that $n$ is nonnegative and that $-\xi\ln \xi\leq \frac{1}{e}$ for all $\xi>0$ we can see that
\begin{align*}
\intomega n|\ln n|&=\F(n,z)-\frac{1}{2}\intomega|\nabla z|^2+\ln\mu\intomega n-2\int_{\{n<1\}}n\ln n\leq \F(n,z)+\ln\mu\intomega n+\frac{2|\Omega|}{e},
\end{align*}
proving \eqref{eq:Fineq_n}. Similarly, we may compute
\begin{align*}
\frac{1}{2}\intomega|\nabla z|^2=\F(n,z)-\mu\intomega\frac{n}{\mu}\ln\frac{n}{\mu}\leq\F(n,z)+\frac{\mu|\Omega|}{e},
\end{align*} 
which first proves \eqref{eq:Fineq_z} and, upon reordering and dropping the nonnegative term, also \eqref{eq:Fineq_below}.
\end{bew}

The main ingredient in showing that this generalized energy is non-increasing (after some waiting time) will be the following differential inequality.

\begin{lemma}\label{lem:F-diffineq}
Let $m>0$ and $T\geq0$ and assume that for $f\in C^{3}\!\left([0,\infty)\right)$ satisfying \eqref{eq:fnoepsprop} the triple $(n,z,u)\in\CSp{2,1}{\bomega\times(T,\infty)}$ is a classical solution of \eqref{CSznoeps}--\,\eqref{CSznoepsBC} in $\Omega\times(T,\infty)$ satisfying $\intomega|u(\cdot,T)|^4\leq \ell$, and $\intomega n(\cdot,t)\leq m$ for all $t>T$, as well as $n>0$ in $\Omega\times(T,\infty)$. Then for all $\mu>0$ we have
\begin{align*}
\frac{\intd}{\intd t}\F\big(n(\cdot,t),z(\cdot,t)\big)&+\intomega\frac{|\nabla n(\cdot,t)|^2}{n(\cdot,t)}\\&+\left\{\frac{1}{2}-\frac{K_3}{2}\intomega|\nabla z(\cdot,t)|^2-K_3^2K_{u}|\Omega|^\frac{1}{4}\big(\ell e^{-\lambda_1 (t-T)}+m\big)\right\}\intomega|\Delta z(\cdot,t)|^2\leq 0
\end{align*}
for all $t>T$, with $K_3$ as in \eqref{global_gnb_constants3} and $K_{u}$, $\lambda_1$ provided by Lemma \ref{lem:l4_bound_u}.
\end{lemma}

\begin{bew}
Since $n$ is positive in $\bomega\times(T,\infty)$ we see by utilizing integration by parts that
\begin{align}\label{eq:F-diffineq1}
\frac{\intd}{\intd t}\F(n,z)&
=-\intomega\frac{|\nabla n|^2}{n}-\intomega|\Delta z|^2+\intomega\Delta z|\nabla z|^2-\intomega\Delta z(u\cdot \nabla z)
\end{align}
holds for all $t>T$, where we used the first and second equations of \eqref{CSznoeps} and $\dive u=0$. By Young's inequality and \eqref{global_gnb_constants3} we have
\begin{align}\label{eq:F-diffineq2}
\intomega \Delta z|\nabla z|^2&\leq\frac{1}{2}\intomega|\Delta z|^2+\frac{1}{2}\intomega|\nabla z|^4
\leq \left\{\frac{1}{2}+\frac{K_3}{2}\intomega|\nabla z|^2\right\}\intomega|\Delta z|^2\quad\text{for all }t>T.
\end{align}
To estimate the last term in \eqref{eq:F-diffineq1}, we note that by Hölder's inequality and \eqref{global_gnb_constants3} there holds $\|\nabla z\|_{\Lo[4]}\leq K_3^2|\Omega|^{\frac{1}{4}}\|\Delta z\|_{\Lo[2]}$ for all $t>T$, which together with Lemma \ref{lem:l4_bound_u} implies
\begin{align}\label{eq:F-diffineq3}
\intomega |\Delta z(u\cdot\nabla z)|&\leq\|\Delta z\|_{\Lo[2]}\|u\|_{\Lo[4]}\|\nabla z\|_{\Lo[4]}\nonumber\\
&\leq K_3^2|\Omega|^\frac{1}{4}\|\Delta z\|_{\Lo[2]}^2\|u\|_{\Lo[4]}\nonumber\\
&\leq K_3^2K_{u}|\Omega|^\frac{1}{4}\big(\ell e^{-\lambda_1 (t-T)}+m\big)\intomega|\Delta z|^2\quad\text{for all }t>T,
\end{align}
since $\intomega n\leq m$ in $(T,\infty)$. Combining \eqref{eq:F-diffineq1}--\eqref{eq:F-diffineq3} and reordering appropriately completes the proof.
\end{bew}

In view of the lemma above, the possibility for an inequality of the form $\frac{\intd}{\intd t}\F\big(n(\cdot,t),z(\cdot,)\big)\leq 0$ will depend on the nonnegativity of the term $\frac{1}{2}-\frac{K_3}{2}\intomega|\nabla z(\cdot,t)|^2-K_3^2K_u|\Omega|^\frac{1}{4}(\ell e^{-\lambda_1t}+m)$. Most of all, this will require some large waiting time $t_0$ and some small bound on $\intomega n$ in order to treat the term $\ell e^{-\lambda_1 (t-T)}+m$. Similarly to the fluid free case, we further require that the energy at a certain time is already sufficiently small, which will provide control of the term containing $\intomega |\nabla z|^2$.

\begin{lemma}\label{lem:decreasing-energy}
Let $T\geq0$ and $\left(4K_3^2K_{u}|\Omega|^\frac{1}{4}\right)^{-1}>m_0>0$, with $K_3$ and $K_{u}$ provided by \eqref{global_gnb_constants3} and Lemma \ref{lem:l4_bound_u}, respectively. Suppose that for $f\in C^{3}\!\left([0,\infty)\right)$ satisfying \eqref{eq:fnoepsprop} the triple $(n,z,u)\in\CSp{2,1}{\bomega\times(T,\infty)}$ is a classical solution of \eqref{CSznoeps}--\,\eqref{CSznoepsBC} in $\Omega\times(T,\infty)$ satisfying $\intomega |u(\cdot,T)|^4\leq\ell$ and $m:=\intomega n(\cdot,T)\leq m_0$, as well as $n>0$ in $\Omega\times(T,\infty)$ and $z\in\CSp{0}{[T,\infty);\W[1,2]}$. Then if
there exist $t_0\geq T$ and $\mu>0$ such that
\begin{align}\label{eq:dec-energy-t0}
\ell e^{-\lambda_1 (t_0-T)}+m_0\leq\frac{1}{4K_3^2K_{u}|\Omega|^\frac{1}{4}}
\end{align}
and
\begin{align}\label{eq:dec-energy-F-bound}
\F\big(n(\cdot,t_0),z(\cdot,t_0)\big)<\frac{1}{4K_3}-\frac{\mu|\Omega|}{e},
\end{align}
then
\begin{align}\label{eq:dec-energy-F}
\frac{\intd}{\intd t}\F\big(n(\cdot,t),z(\cdot,t)\big)\leq0\quad\text{for all }t>t_0.
\end{align}
Furthermore, one can find $\kappa>0$ such that 
\begin{align}\label{eq:ev-spatiotempbounds}
\int_{t_0}^t\!\intomega\frac{|\nabla n|^2}{n}+\kappa\int_{t_0}^t\!\intomega|\Delta z|^2<\frac{1}{4K_3}\quad\text{for all }t>t_0.
\end{align}
\end{lemma}

\begin{bew}
First we note that in view of Remark \ref{rem:masscons-sys-csz} the inequality in \eqref{eq:dec-energy-t0} implies that
\begin{align}\label{eq:t0-est}
\ell e^{-\lambda_1 (t-T)}+m< \ell e^{-\lambda_1 (t_0-T)}+m_0\leq\frac{1}{4K_3^2K_{u}|\Omega|^{\frac{1}{4}}}\quad\text{for all }t>t_0.
\end{align}
Furthermore, recalling Lemma \ref{lem:Fprop} we see that \eqref{eq:dec-energy-F-bound} implies $\frac{K_3}{2}\intomega|\nabla z(\cdot,t_0)|^2\leq K_3\F\big(n(\cdot,t_0),z(\cdot,t_0)\big)+\frac{K_3\mu|\Omega|}{e}<\frac{1}{4}$. Therefore, the set
\begin{align*}
S:=\Big\{T'>t_0\,\Big\vert\, \frac{K_3}{2}\intomega|\nabla z(\cdot,t)|^2<\frac{1}{4}\ \text{ for all }t\in[t_0,T')\Big\}
\end{align*}
is not empty and $T_S:=\sup S$ is a well-defined element of $(t_0,\infty]$. In order to verify that actually $T_S=\infty$ we assume $T_S<\infty$ an derive a contradiction. To this end, we make use of Lemma \ref{lem:F-diffineq} to obtain from the definition of $T_S$ and \eqref{eq:t0-est} that
\begin{align}\label{eq:F-diffineq-kappa}
\frac{\intd}{\intd t}\F\big(n(\cdot,t),z(\cdot,t)\big)&+\intomega\frac{|\nabla n(\cdot,t)|^2}{n(\cdot,t)}+\kappa\intomega|\Delta z(\cdot,t)|^2\leq 0\quad\text{for all }t\in(t_0,T_S),
\end{align}
with some small $\kappa>0$. Due to the assumed $\W[1,2]$-valued continuity of $z$, the mapping $[t_0,\infty)\ni t\mapsto\F\big(n(\cdot,t),z(\cdot,t)\big)$ is continuous as well and we infer from the definition of $T_S$ that $\frac{K_3}{2}\intomega|\nabla z|^2<\frac{1}{4}$ for all $t\in(t_0,T_S)$, but 
\begin{align}\label{eq:F-contradict}
\frac{K_3}{2}\intomega|\nabla z(\cdot,T_S)|^2=\frac{1}{4}.
\end{align}
Integrating \eqref{eq:F-diffineq-kappa} we obtain
\begin{align*}
\F\big(n(\cdot,T_S),z(\cdot,T_S))\leq \F\big(n(\cdot,t_0),z(\cdot,t_0)\big),
\end{align*}
which by Lemma \ref{lem:Fprop} and \eqref{eq:dec-energy-F-bound} shows
\begin{align*}
\intomega|\nabla z(\cdot,T_S)|^2&\leq 2\F\big(n(\cdot,T_S),z(\cdot,T_S))+\frac{2\mu|\Omega|}{e}\leq2\F\big(n(\cdot,t_0),z(\cdot,t_0)\big)+\frac{2\mu|\Omega|}{e}<\frac{1}{2K_3},
\end{align*}
contradicting \eqref{eq:F-contradict} and thus proving $T_S=\infty$. Therefore, the inequality \eqref{eq:F-diffineq-kappa} actually holds for all $t>t_0$, which firstly proves \eqref{eq:dec-energy-F} and secondly, upon integration of \eqref{eq:F-diffineq-kappa} shows \eqref{eq:ev-spatiotempbounds} due to \eqref{eq:dec-energy-F-bound}. 
\end{bew}

\subsection{Conditional regularity estimates}\label{sec4:regularity}
In this section we will establish appropriate Hölder bounds for the components of our approximate solutions under the assumption that we already have control of $\intomega|\nabla z|^p$ for some $p>2$. In fact, as we will see in Section \ref{sec4:obtainingbound}, obtaining the bound assumed throughout the section for the special value of $p=4$, will only require bounds on $\intomega n|\ln n|$ and $\intomega|\nabla z|^2$, which (at least for possibly large times) can be obtained by relying on our analysis of $\F$ (see Section \ref{sec4:evs}). Our arguments here are inspired an approach illustrated in \cite[Section 4.2 and 4.3]{Win16CS2}.
\begin{lemma}\label{lem:cond-reg-n-inf}
Let $p>2$, $m_0>0$, $M>0$ and $\tau>0$. Then there exists $C=C(p,m_0,M,\tau)>0$ such that if for $f\in C^3\!\left([0,\infty)\right)$ satisfying \eqref{eq:fnoepsprop} and some $t_0\geq0$ the triple $(n,z,u)\in\CSp{2,1}{\bomega\times(t_0,\infty)}$ is a classical solution of \eqref{CSznoeps}--\eqref{CSznoepsBC} in $\Omega\times(t_0,\infty)$ satisfying $n\geq0$ in $\Omega\times(t_0,\infty)$ and
\begin{align}\label{eq:cond-reg-n-inf-nl1}
\intomega n(\cdot,t_0)\leq m_0
\end{align}
as well as
\begin{align*}
\intomega |\nabla z(\cdot,t)|^p\leq M\quad\text{for all }t>t_0,
\end{align*}
then
\begin{align}\label{eq:cond-reg-n-inf-ninf}
\|n(\cdot,t)\|_{\Lo[\infty]}\leq C\quad\text{for all }t\geq t_0+\tau.
\end{align}
\end{lemma}

\begin{bew}
The proof is based on arguments employed in e.g. \cite[Lemma 4.4]{Win16CS2}. We let $T>t_0+1$ and define
\begin{align*}
S(T):=\max\left\{S_1,S_2(T)\right\}
\end{align*}
with
\begin{align*}
S_1:=\max_{t\in[t_0,t_0+1]}(t-t_0)\|n(\cdot,t)\|_{\Lo[\infty]}\quad\text{and}\quad S_2(T):=\max_{t\in[t_0+1,T]}\|n(\cdot,t)\|_{\Lo[\infty]}.
\end{align*}
Now, in order to estimate $S(T)$ from above, we let $t_1(t):=\max\{t-1,t_0\}$ and for $t\in(t_0,T)$ represent $n(\cdot,t)$ according to
\begin{align}\label{eq:cond-ninf-varconstants}
n(\cdot,t)&=e^{(t-t_1)\Delta}n(\cdot,t_1)+\!\int_{t_1}^t\!\!e^{(t-s)\Delta}\Big[\nabla\divedot\big(n(\cdot,s)f'(n(\cdot,s))\nabla z(\cdot,s)\big)-\big(u(\cdot,s)\divedot\nabla n(\cdot,s)\big)\Big]\intd s\nonumber\\
&=:e^{(t-t_1)\Delta}n(\cdot,t_1)+I(t_1,t),
\end{align}
where $(e^{\sigma\Delta})_{\sigma\geq0}$ denotes the heat semigroup with Neumann boundary data in $\Omega$. Fixing some $q\in(2,p)$, we may rely on well known estimates for the heat semigroup (e.g. \cite[Lemma 1.3]{win10jde} and \cite[Lemma 3.3]{FIYW16}) to find $C_1>0$ and $C_2>0$ such that for all $\sigma\in(0,1)$ there holds
\begin{align}\label{eq:cond-ninf-heatsemi1}
\|e^{\sigma\Delta}\varphi\|_{\Lo[\infty]}\leq C_1\sigma^{-1}\|\varphi\|_{\Lo[1]}\quad\text{for all }\varphi\in\Lo[1]
\end{align}
and
\begin{align}\label{eq:cond-ninf-heatsemi2}
\|e^{\sigma\Delta}\nabla\cdot\varphi\|_{\Lo[\infty]}\leq C_2\sigma^{-\gamma}\|\varphi\|_{\Lo[q]}\quad\text{for all }\varphi\in \CSp{1}{\bomega}\text{ such that }\varphi\cdot\nu=0\text{ on }\romega,
\end{align}
with $\gamma:=\frac{1}{2}+\frac{1}{q}<1$. In the case $t\in(t_0,t_0+1]$, when $t_1(t)=t_0$, we thus have
\begin{align}\label{eq:cond-ninf-ineq0}
\big\|e^{(t-t_0)\Delta}n(\cdot,t_0)\big\|_{\Lo[\infty]}\leq C_1 m_0(t-t_0)^{-1},
\end{align}
thanks to \eqref{eq:cond-reg-n-inf-nl1} and \eqref{eq:cond-ninf-heatsemi1}. Furthermore, making use of $\dive u=0$, the fact that $f'\leq 1$ on $[0,\infty)$, and \eqref{eq:cond-ninf-heatsemi2} we see that
\begin{align*}
&\|I(t_0,t)\|_{\Lo[\infty]}\leq C_2\int_{t_0}^t(t-s)^{-\gamma}\big(\big\|n(\cdot,s)\nabla z(\cdot,s)\big\|_{\Lo[q]}+\big\|n(\cdot,s)u(\cdot,s)\big\|_{\Lo[q]}\big)\intd s
\end{align*}
holds for all $t\in(t_0,t_0+1]$. Herein, multiple applications of the Hölder inequality show that
\begin{align}\label{eq:cond-ninf-hoelder1}
\nonumber\big\|n(\cdot,s)\nabla z(\cdot,s)\big\|_{\Lo[q]}
&\leq\|n(\cdot,s)\|_{\Lo[\infty]}^a\|n(\cdot,s)\|_{\Lo[1]}^{1-a}\|\nabla z(\cdot,s)\|_{\Lo[p]}\\&\leq m_0^{1-a}M^{\frac{1}{p}}\|n(\cdot,s)\|_{\Lo[\infty]}^a\quad\text{for all }s>t_0
\end{align}
with $a:=1-\frac{p-q}{pq}\in(0,1)$ and 
\begin{align}\label{eq:cond-ninf-hoelder2}
\big\|n(\cdot,s)u(\cdot,s)\big\|_{\Lo[q]}\leq C_3(1+m_0) m_0^{1-a}\|n(\cdot,s)\|_{\Lo[\infty]}^a\quad\text{for all }s>t_0,
\end{align}
for some $C_3>0$, where $\|u(\cdot,t)\|_{\Lo[p]}\leq C_3(1+m_0)$ in view of Lemma \ref{lem:stokes_est}. In particular, recalling the definition of $S_1$ we have
\begin{align}\label{eq:cond-ninf-ineq1}
\|I(t_0,t)\|_{\Lo[\infty]}\leq C_4 S_1^a\int_{t_0}^t(t-s)^{-\gamma}(s-t_0)^{-a}\intd s\quad\text{for all }t\in(t_0,t_0+1].
\end{align}
with some $C_4>0$. Since $\int_{t_0}^t(t-s)^{-\gamma}(s-t_0)^{-a}\intd s=(t-t_0)^{1-\gamma-a}\int_0^1(1-\zeta)^{-\gamma}\zeta^{-a}\intd \zeta\leq B(1-a,1-\gamma)$ is finite according to the facts that $a<1$ and $\gamma<1$, we consequently see that collecting \eqref{eq:cond-ninf-varconstants}, \eqref{eq:cond-ninf-ineq0}, and \eqref{eq:cond-ninf-ineq1} shows that there exists some $C_5>0$ such that
\begin{align*}
(t-t_0)\|n(\cdot,t)\|_{\Lo[\infty]}\leq C_5+C_5 S_1^a\quad\text{for all }t\in(t_0,t_0+1],
\end{align*}
which, due to $a<1$, implies that
\begin{align}\label{eq:cond-ninf-S1bound}
S_1\leq C_6:=\max\big\{1,(2C_5)^{\frac{1}{1-a}}\big\}.
\end{align}
The estimation of $S_2(T)$ follows a similar path. We fix $t\in[t_0+1,T]$ and obtain from \eqref{eq:cond-ninf-varconstants}, \eqref{eq:cond-ninf-heatsemi1}, and \eqref{eq:cond-ninf-heatsemi2} that
\begin{align*}
\|n(\cdot,t)\|_{\Lo[\infty]}&\leq \big\|e^{\Delta}n(\cdot,t-1)\big\|_{\Lo[\infty]}+\|I(t-1,t)\|_{\Lo[\infty]}\\
&\leq C_1\|n(\cdot,t-1)\|_{\Lo[1]}+C_2\int_{t-1}^t(t-s)^{-\gamma}\big(\big\|n(\cdot,s)\nabla z(\cdot,s)-n(\cdot,s)u(\cdot,s)\big\|_{\Lo[q]}\big)\intd s.
\end{align*}
From which, again by relying on \eqref{eq:cond-reg-n-inf-nl1}, \eqref{eq:cond-ninf-hoelder1}, and \eqref{eq:cond-ninf-hoelder2}, we infer that
\begin{align*}
\|n(\cdot,t)\|_{\Lo[\infty]}\leq C_1m_0+C_2m_0^{1-a}\big(M^{\frac{1}{p}}+C_3(1+m_0)\big)\int_{t-1}^t(t-s)^{-\gamma}\|n_\eps(\cdot,s)\|_{\Lo[\infty]}^a\intd s	
\end{align*}
holds for all $t\in[t_0+1,T]$. By the definition of $S_2(T)$ we have $\|n(\cdot,s)\|_{\Lo[\infty]}^a\leq S_2^a(T)$ for all $s\in[t_0+1,T]$, so that in both of the cases $t\in[t_0+1,t_0+2]$ and $t>t_0+2$ we may estimate
\begin{align*}
\int_{t-1}^t(t-s)^{-\gamma}\|n(\cdot,s)\|_{\Lo[\infty]}^a\intd s&\leq S_1^a\int_{t-1}^t(t-s)^{-\gamma}(s-t_0)^{-a}\intd s+S_2^a(T)\int_{t-1}^t(t-s)^{-\gamma}\intd s\\&\leq C_7S_{1}^a+\frac{1}{1-\gamma}S_2^a(T).
\end{align*}
with some $C_7>0$. Collecting these estimates and making use of \eqref{eq:cond-ninf-S1bound} we find $C_8>0$ such that
\begin{align*}
\|n(\cdot,t)\|_{\Lo[\infty]}\leq C_8+C_8 S_2^a(T)\quad\text{for all }t\in[t_0+1,T],
\end{align*}
which implies $S_2(T)\leq C_9:=\max\big\{1,(2C_8)^{\frac{1}{1-a}}\big\}$ for all $T>t_0+1$. Finally, combining both estimates for $S_1$ and $S_2(T)$ establishes \eqref{eq:cond-reg-n-inf-ninf} if we let $C:=\max\{S_1,\tfrac{S_1}{\tau},C_9\}$.
\end{bew}

With the improved regularity for $n$ at hand, we can easily derive time local Hölder continuity of $n$ and $u$ under the same assumptions as above.

\begin{lemma}\label{lem:cond-reg-nu-hoelder}
Let $p>2$, $m_0>0$, $M>0$ and $\tau>0$. Then there exist some $\theta=\theta(p)\in(0,1)$ and $C=C(p,m_0,M,\tau)>0$ such that if $f\in C^3\!\left([0,\infty)\right)$ satisfies \eqref{eq:fnoepsprop} and if for some $t_0\geq0$ the triple $(n,z,u)\in\CSp{2,1}{\bomega\times(t_0,\infty)}$ is a classical solution of \eqref{CSznoeps}--\eqref{CSznoepsBC} in $\Omega\times(t_0,\infty)$ with the properties that $n\geq0$ in $\Omega\times(t_0,\infty)$ and
\begin{align}\label{eq:cond-reg-nu-hoelder-nl1}
\intomega n(\cdot,t_0)\leq m_0,
\end{align}
as well as
\begin{align}\label{eq:cond-reg-nu-hoelder-nabzlp}
\intomega |\nabla z(\cdot,t)|^p\leq M\quad\text{for all }t>t_0,
\end{align}
then
\begin{align*}
\|n\|_{\CSp{\theta,\frac{\theta}{2}}{\bomega\times[t,t+1]}}\leq C\quad\text{ and }\quad \|u\|_{\CSp{\theta,\frac{\theta}{2}}{\bomega\times[t,t+1]}}\leq C\quad\text{for all }t\geq t_0+\tau.
\end{align*}
\end{lemma}

\begin{bew}
With $\alpha$ given by \eqref{IR} we fix $\beta\in\big(\frac{1}{2},\alpha\big)$. Then we apply the fractional power $A^\beta$ of the $L^2$--realization of the Stokes operator to a variation-of-constants representation for $u$ to obtain the identity
\begin{align*}
A^\beta u(\cdot,t)=A^\beta e^{-(t-t_1)A}u(\cdot,t_1)+\int_{t_1}^t A^\beta e^{-(t-s)A}\HP\left(n(\cdot,s)\nabla\phi\right)\intd s,\quad t\geq t_1,
\end{align*}
where $t_1:=\max\{t-1,t_0\}$. Recalling that the positive sectorial Stokes operator $A$ generates the contracting semigroup $\big(e^{-tA}\big)_{t\geq0}$ in $L_{\sigma}^2\left(\Omega\right)$ and the fractional powers of the Stokes operator fulfill the decay property
\begin{align*}
\big\|A^\beta e^{-tA}\big\|\leq C_1 t^{-\beta}e^{-\lambda_1 t}\quad\text{for all }t>0,
\end{align*}
with some $C_1>0$ (\cite[Theorem 37.5]{sellyou}), we can make use of the boundedness of $\HP$ in $\Lo[2]$, \eqref{phireg}, \eqref{eq:cond-reg-nu-hoelder-nl1}, and Lemma \ref{lem:stokes_est} to obtain $C_1>0$ such that
\begin{align}\label{eq:con-reg-nu-hoelder-eq1}
\nonumber\big\|A^\beta u(\cdot,t)\big\|_{\Lo[2]}&\leq \big\|A^\beta e^{-(t-t_1)A}u(\cdot,t_1)\big\|_{\Lo[2]}+\int_{t_1}^t\big\|A^\beta e^{-(t-s)A}\HP\left(n(\cdot,s)\nabla\phi\right)\big\|_{\Lo[2]}\intd s\\
&\leq C_1(t-t_1)^{-\beta}+C_1K_1\int_{t_1}^t(t-s)^{-\beta}\|n(\cdot,s)\|_{\Lo[2]}\intd s
\end{align}
for all $t>t_1$. Since the assumptions \eqref{eq:cond-reg-nu-hoelder-nl1} and \eqref{eq:cond-reg-nu-hoelder-nabzlp} allow for an application of Lemma \ref{lem:cond-reg-n-inf}, we can find $C_2>0$ such that $\|n(\cdot,t)\|_{\Lo[2]}\leq C_2$ for all $t\geq t_0+\tau$. Combining $\beta<1$ with the fact that in both cases $(t-t_1)^{1-\beta}\leq1$ and $(t-t_1)^{-\beta}\leq 1+\tau^{-\beta}$ hold for $t\geq t_0+\tau$, we infer from \eqref{eq:con-reg-nu-hoelder-eq1} the existence of some $C_3:=C_3(p,m_0,M,\tau)>0$ such that
\begin{align*}
\big\|A^\beta u(\cdot,t)\big\|_{\Lo[2]}\leq C_3\quad\text{for all }t\geq t_0+\tau.
\end{align*}
Considering that since $\beta\in(\frac{1}{2},\alpha)$ the domains of fractional powers of the Stokes semigroup satisfy $\DA\hookrightarrow D(A^{\beta})\hookrightarrow \CSp{\theta_1}{\bomega}$ for any $\theta_1\in(0,2\beta-1)$ (\cite[Lemma III.2.4.2]{sohr} and \cite[Theorem 5.6.5]{evans}), the previous estimate entails the existence of some $C_4>0$ such that
\begin{align*}
\|u(\cdot,t)\|_{\CSp{\theta_1}{\bomega}}\leq C_4\quad\text{for all }t\geq t_0+\tau.
\end{align*}
Making use of similar arguments we can find $C_5>0$ such that
\begin{align*}
\big\|A^\beta u(\cdot,t)-A^\beta u(\cdot,t_2)\big\|_{\Lo[2]}\leq C_5(t-t_2)^{1-\beta}\quad\text{for all }t_2\geq t_0+\tau\text{ and }t\in[t_2,t_2+1],
\end{align*}
which together with \eqref{eq:con-reg-nu-hoelder-eq1} readily implies the Hölder regularity of $u$ for some $\theta_2:=\min\{1-\beta,\theta_1\}$. For the regularity of $n$ we first note that by Lemma \ref{lem:cond-reg-n-inf} we obtain a constant $C_6:=C_6(p,m_0,M,\tau)$ such that $n(x,t)\leq C_6$ for all $x\in\Omega$ and $t\geq t_0+\frac{\tau}{2}$. Hence, the function $n$ is a bounded distributional solution to the parabolic equation
\begin{align*}
\tilde{n}_{t}-\dive a(x,t,\tilde{n},\nabla \tilde{n})=0\quad \text{in }\Omega\times(t_0,\infty),
\end{align*}
with $a(x,t,\tilde{n},\nabla\tilde{n}):=\nabla\tilde{n}+n(x,t)f'\big(n(x,t)\big)\nabla z(x,t)-u n$ and $a(x,t,\tilde{n},\nabla\tilde{n})\cdot\nu=0$ on the boundary of $\Omega$. Considering that with the arguments illustrated in the first part of the proof we can find $C_7:=C_7(p,m_0,M,\tau)$ such that $|u(x,t)|\leq C_7$ for all $x\in\Omega$ and $t\geq t_0+\frac{\tau}{2}$, we let $\psi_0(x,t):=n(x,t)^2|\nabla z(x,t)|^2+|u(x,t)n(x,t)|^2$ and $\psi_1(x,t):=C_6|\nabla z(x,t)|+C_6C_7$ and then see by means of Young's inequality and \eqref{eq:feps_prop1} that
\begin{align*}
a(x,t,\tilde{n},\nabla\tilde{n})\nabla\tilde{n}\geq\frac{1}{2}|\nabla\tilde{n}|^2-\psi_0
\quad\text{and}\quad
|a(x,t,\tilde{n},\nabla\tilde{n})|\leq|\nabla\tilde{n}(x,t)|+\psi_1(x,t)
\end{align*}
for all $(x,t)\in\Omega\times(t_0+\frac{\tau}{2},\infty)$. Since \eqref{eq:cond-reg-nu-hoelder-nabzlp} provides a bound for $|\nabla z|^2$ in $\LSp{\infty}{(t_0,\infty);\Lo[\frac{p}{2}]}$, we obtain from a well known result in \cite[Theorem 1.3]{PorzVesp93} that $\|n\|_{\CSp{\theta_3,\frac{\theta_3}{2}}{\bomega\times[t,t+1]}}\leq C_8$ for all $t>t_0+\tau$ with some $\theta_3(p)>0$ and $C_8>0$. Picking $\theta\in(0,\min\{\theta_2,\theta_3\})$ the claim follows immediately.
\end{bew}

In order to prepare a further improvement on the regularity we will show the following.

\begin{lemma}\label{lem:cond-reg-z-inf}
Let $p>2$, $m_0>0$, $m_1>0$, $M>0$ and $T>0$. Then there exists $C=C(p,m_0,m_1,M,T)>0$ such that if for $f\in C^3\!\left([0,\infty)\right)$ satisfying \eqref{eq:fnoepsprop} and $t_0\geq0$ the triple $(n,z,u)\in\CSp{0}{\bomega\times[t_0,\infty)}\cap\CSp{2,1}{\bomega\times(t_0,\infty)}$ is a classical solution of \eqref{CSznoeps}--\eqref{CSznoepsBC} in $\Omega\times(t_0,\infty)$ with the properties that $n\geq0$ in $\Omega\times(t_0,\infty)$ and 
\begin{align}\label{eq:cond-reg-z-inf-nl1}
\intomega n(\cdot,t)\leq m_0\quad\text{for all }t>t_0,
\end{align}
and
\begin{align}\label{eq:cond-reg-z-inf-zl1}
\intomega z(\cdot,t_0)\leq m_1,
\end{align}
as well as
\begin{align}\label{eq:cond-reg-z-inf-nabzlp}
\intomega |\nabla z(\cdot,t)|^p\leq M\quad\text{for all }t>t_0,
\end{align}
then
\begin{align*}
z(x,t)\leq C\quad\text{for all }x\in\Omega\ \text{and}\ t\in(t_0,T).
\end{align*}
\end{lemma}

\begin{bew}
Because of the assumption $p>2$ we have $\W[1,p]\hookrightarrow \Co[1-\frac{2}{p}]$ and thus there exists some constant $C_1>0$ such that for each $\varphi\in\W[1,p]$ it holds that
\begin{align}\label{eq:cond-reg-z-inf-hoelder}
|\varphi(x)-\varphi(y)|\leq C_1|x-y|^{1-\frac{2}{p}}\|\nabla \varphi\|_{\Lo[p]}\quad\text{for all }x,y\in\Omega.
\end{align}
By Lemma \ref{lem:bound_zeps}, Remark \ref{rem:masscons-sys-csz} and the assumptions \eqref{eq:cond-reg-z-inf-nl1} and \eqref{eq:cond-reg-z-inf-zl1} we see that
\begin{align*} 
\intomega z(\cdot,t)\leq \intomega z(\cdot,t_0)+m_0(t-t_0)\leq m_1+m_0T\quad\text{for all }t\in(t_0,T),
\end{align*}
whence for any such $t\in(t_0,T)$ we can find $x_0(t)\in\Omega$ such that
\begin{align*}
z(x_0(t),t)\leq\frac{m_1+m_0T}{|\Omega|}.
\end{align*}
Therefore, \eqref{eq:cond-reg-z-inf-hoelder} in conjunction with the assumption \eqref{eq:cond-reg-z-inf-nabzlp} shows that
\begin{align*}
z(x,t)&\leq z(x_0(t),t)+\big|z(x,t)-z(x_0(t),t)\big|\\
&\leq \frac{m_1+m_0T}{|\Omega|}+C_1|x-x_0(t)|^{1-\frac{2}{p}}\|\nabla z(\cdot,t)\|_{\Lo[p]}\\
&\leq \frac{m_1+m_0T}{|\Omega|}+C_2 M^\frac{1}{p}
\end{align*}
holds for all $x\in\Omega$, with $C_2$ only depending on $p$ and the diameter of $\Omega$.
\end{bew}

Drawing on the now proven time-local bound for $z$, we can rely on the Hölder estimates for $n$ and $u$ and well known parabolic regularity theory to the following set of further bounds.

\begin{lemma}\label{lem:cond-reg-hoelder2}
Let $p>2,m_0>0,m_1>0,M>0,T>0$ and $\tau>0$. Then there exist $\theta=\theta(p)\in(0,1)$ and $C=C(p,m_0,m_1,M,T,\tau)>0$ such that if for $f\in C^3\!\left([0,\infty)\right)$ satisfying \eqref{eq:fnoepsprop} and $t_0\geq0$ the triple $(n,z,u)\in\CSp{0}{\bomega\times[t_0,\infty)}\cap\CSp{2,1}{\bomega\times(t_0,\infty)}$ is a classical solution of \eqref{CSznoeps}--\eqref{CSznoepsBC} in $\Omega\times(t_0,\infty)$ with the properties that $n\geq0$ and $z\geq0$ in $\Omega\times(t_0,\infty)$ and
\begin{align*}
\intomega n(\cdot,t_0)\leq m_0,
\end{align*}
and
\begin{align*}
\intomega z(\cdot,t_0)\leq m_1,
\end{align*}
as well as
\begin{align*}
\intomega |\nabla z(\cdot,t)|^p\leq M\quad\text{for all }t>t_0,
\end{align*}
then
\begin{align}\label{eq:cond-reg-hoelder2-bounds}
\|n\|_{\CSp{2+\theta,1+\frac{\theta}{2}}{\bomega\times[t_0+\tau,T]}}\leq C,\quad\|z\|_{\CSp{2+\theta,1+\frac{\theta}{2}}{\bomega\times[t_0+\tau,T]}}\leq C,\quad \|u\|_{\CSp{2+\theta,1+\frac{\theta}{2}}{\bomega\times[t_0+\tau,T]}}\leq C.
\end{align}
\end{lemma}

\begin{bew}
By Lemma \ref{lem:cond-reg-z-inf} and the fact that $z$ is nonnegative we have
\begin{align*}
0\leq z\leq C_1\quad\text{in }\Omega\times(t_0,T)
\end{align*}
with some $C_1=C_1(p,m_0,m_1,M,T)>0$. Thus, letting $\tilde{c}:=e^{-z}$ we obtain
\begin{align}\label{eq:cond-reg-hoelder2-tildezeps-bound}
e^{-C_1}\leq\tilde{c}\leq 1\quad\text{in }\Omega\times(t_0,T).
\end{align}
Furthermore, $\tilde{c}$ solves the Neumann boundary value problem $\tilde{c}_t=\Delta \tilde{c}+u\nabla\tilde{c}-f(n)\tilde{c}$ in $\Omega\times(t_0,\infty)$ with Hölder continuous coefficients, since Lemma \ref{lem:cond-reg-nu-hoelder} entails the existence of $\theta_1\in(0,1)$ and $C_2=C_2(p,m_0,M,\tau)>0$ such that
\begin{align*}
\|n\|_{\CSp{\theta_1,\frac{\theta_1}{2}}{\bomega\times[t_0+\frac{\tau}{4},T]}}+\|u\|_{\CSp{\theta_1,\frac{\theta_1}{2}}{\bomega\times[t_0+\frac{\tau}{4},T]}}\leq C_2.
\end{align*}
Hence, according to standard parabolic Schauder theory (\cite[III.5.1 and IV.5.3]{LSU}), there exists some $\theta_2\in(0,1)$ and $C_3=C_3(p,m_0,m_1,M,T,\tau)$ such that
\begin{align*}
\|\tilde{c}\|_{\CSp{2+\theta_2,1+\frac{\theta_2}{2}}{\bomega\times[t_0+\frac{\tau}{2},T]}}\leq C_3,
\end{align*}
yielding the regularity assertion for $z$ featured in \eqref{eq:cond-reg-hoelder2-bounds} due to the lower bound for $\tilde{c}$ in \eqref{eq:cond-reg-hoelder2-tildezeps-bound}. Relying on parabolic Schauder theory once more, we can conclude from the first equation that also $n$ satisfies \eqref{eq:cond-reg-hoelder2-bounds}. That also $u$ satisfies \eqref{eq:cond-reg-hoelder2-bounds} can be readily obtained by well known smoothing properties of the Stokes operator (see eg. \cite[Theorem 2.8]{GigSohr91}, \cite[Theorem 1.1]{Amann00}) and the boundedness of $n$ established in Lemma \ref{lem:cond-reg-n-inf}. 
\end{bew}

\subsection{Conditional estimates for \texorpdfstring{$\intomega|\nabla z|^4$}{the L4-norm of gradient z} and \texorpdfstring{$\intomega n^2$}{the L2-norm of n}}\label{sec4:obtainingbound}
In this section we will focus on obtaining a bound on $\intomega|\nabla z|^4$, which in view of Section \ref{sec4:regularity} is the main requirement for the regularity estimates we will depend on later. As a preliminary step we derive some basic differential inequalities through standard testing procedures.

\begin{lemma}\label{lem:neps_ineq-l2}
Suppose that for $f\in C^3\!\left([0,\infty)\right)$ satisfying \eqref{eq:fnoepsprop} and $t_0\geq0$ the triple $(n,z,u)\in\CSp{2,1}{\bomega\times(t_0,\infty)}$ is a classical solution of \eqref{CSznoeps}--\,\eqref{CSznoepsBC} in $\Omega\times(t_0,\infty)$. Then
\begin{align}\label{eq:neps_ineq-l2}
\frac{\intd}{\intd t}\intomega n^2+\intomega|\nabla n|^2\leq\intomega n^2|\nabla z|^2\quad\text{for all }t>t_0.
\end{align}
\end{lemma}

\begin{bew}
By simply testing the first equation of \eqref{CSznoeps} with $n$, we can rely on integration by parts, one application of Young's inequality, and the fact $|f'(n)|\leq 1$ to easily arrive at \eqref{eq:neps_ineq-l2}.
\end{bew}

\begin{lemma}\label{lem:nab_zeps_ineq-l4}
For any $\eta\in(0,\frac{5}{4})$ there exists $C>0$ such that if for $f\in C^3\!\left([0,\infty)\right)$ satisfying \eqref{eq:fnoepsprop} and $t_0\geq0$ the triple $(n,z,u)\in\CSp{2,1}{\bomega\times(t_0,\infty)}$ is a classical solution of \eqref{CSznoeps}--\,\eqref{CSznoepsBC} in $\Omega\times(t_0,\infty)$ with $n\geq0$ in $\Omega\times(t_0,\infty)$, then 
\begin{align} \label{eq:nab_zeps_ineq-l4}
\frac{\intd}{\intd t}\intomega|\nabla z|^4&+\left(\frac{5}{2}-2\eta\right)\intomega\Big|\nabla|\nabla z|^2\Big|^2\nonumber
\\&\leq 8\intomega|\nabla z|^6+\frac{12}{\eta}\intomega n^2|\nabla z|^2+4\intomega |\nabla z|^4|\nabla u|+C\left(\intomega|\nabla z|^2\right)^2
\end{align}
holds for all $t>t_0$.
\end{lemma}

\begin{bew}
We differentiate the second equation of \eqref{CSznoeps} with regard to space and multiply by $|\nabla z|^2\nabla z$. In the resulting equality we can employ the identity $\nabla z\cdot\nabla\Delta z=\frac{1}{2}\Delta|\nabla z|^2-|D^2z|^2$ to obtain upon integration by parts that
\begin{align}\label{eq:nab_zeps_ineq1}
\nonumber\frac{\intd}{\intd t}\intomega|\nabla z|^4&+2\intomega\big|\nabla|\nabla z|^2\big|^2+4\intomega|\nabla z|^2|D^2z|^2\\&=-4\intomega|\nabla z|^2\nabla z\cdot\nabla|\nabla z|^2-4\intomega|\nabla z|^2f(n)\Delta z-4\intomega f(n)\nabla|\nabla z|^2\cdot\nabla z\nonumber\\&\quad\; -4\intomega|\nabla z|^2\nabla z\cdot(\nabla u\cdot\nabla z)+2\intromega|\nabla z|^2\frac{\partial|\nabla z|^2}{\partial\nu}
\end{align}
holds for all $t>t_0$, due to the fact that $u$ is divergence free and the assumed boundary conditions. Relying on the facts that  $\frac{\partial|\nabla z|^2}{\partial\nu}\leq C_1|\nabla z|^2$ on $\romega$ holds for some $C_1>0$ only depending on $\Omega$ (\cite[Lemma 4.2]{MS14}) and that for fixed $\eta\in(0,\frac{5}{4})$ there exists $C_2>0$ such that $\||\nabla z|^2\|_{\LSp{2}{\romega}}\leq\eta\|\nabla|\nabla z|^2\|_{\Lo[2]}+C_2\|\nabla z\|_{\Lo[2]}$ (c.f. \cite[Remark 52.9]{QS07}), we obtain
\begin{align}\label{eq:nab_zeps_ineq2}
2\intromega|\nabla z|^2\frac{\partial|\nabla z|^2}{\partial\nu}\leq\eta\intomega\big|\nabla|\nabla z|^2\big|^2+C_3\Big(\intomega|\nabla z|^2\Big)^2\quad\text{for all }t>t_0,
\end{align}
with some $C_3>0$. For the remaining integrals, we note that since $f(n)\leq n$ and $|\Delta z|^2\leq 2|D^2z|^2$ by the Cauchy-Schwarz inequality, we can employ Young's inequality to see that
\begin{align}\label{eq:nab_zeps_ineq3}
-4\intomega|\nabla z|^2\nabla z\cdot\nabla|\nabla z|^2&\leq\frac{1}{2}\intomega\big|\nabla|\nabla z|^2\big|^2+8\intomega|\nabla z|^6\quad\text{for all }t>t_0,\\\label{eq:nab_zeps_ineq4}
-4\intomega|\nabla z|^2f(n)\Delta z&\leq\eta\intomega|\nabla z|^2|\Delta z|^2+\frac{4}{\eta}\intomega n^2|\nabla z|^2\nonumber\\&\leq2\eta\intomega|\nabla z|^2|D^2 z|^2+\frac{4}{\eta}\intomega n^2|\nabla z|^2\quad\text{for all }t>t_0,
\intertext{as well as}\label{eq:nab_zeps_ineq5}
-4\intomega f(n)\nabla|\nabla z|^2\cdot\nabla z&\leq\frac{\eta}{2}\intomega\big|\nabla|\nabla z|^2\big|^2+\frac{8}{\eta}\intomega n^2|\nabla z|^2\quad\text{for all }t>t_0.
\end{align}
Collecting \eqref{eq:nab_zeps_ineq1}--\eqref{eq:nab_zeps_ineq5} we thus obtain
\begin{align*}
\frac{\intd}{\intd t}\intomega|\nabla z|^4&+\left(\frac{3}{2}-\frac{3}{2}\eta\right)\intomega\big|\nabla|\nabla z|^2\big|^2+\left(4-2\eta\right)\intomega|\nabla z|^2|D^2z|^2\\&\leq8\intomega|\nabla z|^6+\frac{12}{\eta}\intomega n^2|\nabla z|^2+4\intomega|\nabla z|^4|\nabla u|+C_3\Big(\intomega|\nabla z|^2\Big)^2\quad\text{for all }t>t_0.
\end{align*}
Due to the pointwise inequality $\big|\nabla|\nabla z|^2\big|^2\leq 4|D^2 z|^2|\nabla z|^2$ this readily implies \eqref{eq:nab_zeps_ineq-l4}.
\end{bew}

Combination of the two prepared inequalities will now result in the desired bounds for $\intomega|\nabla z|^4$ and $\intomega n^2$, if we assume that we already have suitable bounds for the quantities $\intomega n\ln n$ and $\intomega|\nabla z|^2$. The bounds on these quantities will later on be obtained from the energy functional upon the requirement that $\intomega n_0$ is small.

\begin{lemma}\label{lem:cond_n-l2_nabz-l4_bound}
Let $K_2$ be as in \eqref{global_gnb_constants1}. Then for all $m_0>0$, each $L>0$ and any $M\in\big(0,\frac{1}{4K_2}\big)$ and $\tau>0$ there exists $C>0$ such that if for $f\in C^3\!\left([0,\infty)\right)$ satisfying \eqref{eq:fnoepsprop} and some $t_0\geq0$ the triple $(n,z,u)\in\CSp{2,1}{\bomega\times(t_0,\infty)}$ is a classical solution of \eqref{CSznoeps}--\,\eqref{CSznoepsBC} in $\Omega\times(t_0,\infty)$ satisfying $n\geq0$ in $\Omega\times(t_0,\infty)$ and
\begin{align}\label{eq:cond_n-l2_m-req}
\intomega n(\cdot,t_0)\leq m_0,
\end{align}
as well as 
\begin{align}\label{eq:cond_n-l2_nlnn-req}
\intomega n(\cdot,t)|\ln n(\cdot,t)|\leq L\quad\text{and}\quad\intomega|\nabla z(\cdot,t)|^2\leq M\quad\text{for all }t>t_0,
\end{align}
then
\begin{align}\label{eq:cond_n-l2_n-l2}
\intomega n^2(\cdot,t)\leq C\quad\text{and}\quad\intomega |\nabla z(\cdot,t)|^4\leq C\quad\text{for all }t\geq t_0+\tau.
\end{align}
\end{lemma}

\begin{bew}
First, we note that due to $M<\frac{1}{4K_2}$, by continuity, one can find some small $\eta\in(0,1)$ such that
\begin{align}\label{eq:cond_n-l2_Mbound}
M<\frac{(2-2\eta)(1-\eta)}{8K_2(1+\eta)}.
\end{align}
Now, assuming \eqref{eq:cond_n-l2_m-req} and \eqref{eq:cond_n-l2_nlnn-req} to hold, we combine the inequalites established in Lemma \ref{lem:neps_ineq-l2} and Lemma \ref{lem:nab_zeps_ineq-l4} to obtain
\begin{align}\label{eq:cond_n-l2-combineq}
\quad&\frac{\intd}{\intd t}\left\{\intomega n^2+\intomega|\nabla z|^4\right\}+\intomega|\nabla n|^2+\Big(\frac{5}{2}-2\eta\Big)\intomega\big|\nabla|\nabla z|^2\big|^2\\\leq\quad &\Big(1+\frac{12}{\eta}\Big)\intomega n^2|\nabla z|^2+8\intomega|\nabla z|^6+4\intomega|\nabla z|^4|\nabla u|+C_1M^2\quad\text{for all }t>t_0,\nonumber
\end{align}
with some $C_1>0$. Herein, Young's inequality provides $C_2>0$ such that 
\begin{align}\label{eq:cond_n-l2-eq1}
\left(1+\frac{12}{\eta}\right)\intomega n^2|\nabla z|^2\leq 8\eta\intomega|\nabla z|^6+C_2\intomega n^3\quad\text{for all }t>t_0.
\end{align}
To further control the term containing $n^3$, we recall that by a variant of the \GNI\ (c.f. \cite[(22)]{Bil94}) and Remark \ref{rem:masscons-sys-csz} we have
\begin{align}\label{eq:cond_n-l2-eq2}
C_2\intomega n^3&\leq\frac{1}{2L}\left(\intomega|\nabla n|^2\right)\left(\intomega n|\ln n|\right)+C_3\left(\intomega n\right)^3+C_3\nonumber\\
&\leq \frac{1}{2}\intomega|\nabla n|^2+C_3m_0^3+C_3\quad\text{for all }t>t_0,
\end{align}
with some $C_3>0$. Returning to the analyzation of the remaining terms in \eqref{eq:cond_n-l2-combineq}, we observe that by Hölder's inequality, Lemma \ref{lem:stokes_est} combined with \eqref{eq:cond_n-l2_m-req}, the \GNI,\ and finally Young's inequality we can find $C_4,C_5,C_6>0$ such that
\begin{align}\label{eq:cond_n-l2-eq3}
4\intomega|\nabla z|^4|\nabla u|&\leq 4\big\||\nabla z|^2\big\|_{\Lo[6]}^2\|\nabla u\|_{\Lo[\frac{3}{2}]}\leq C_4(1+m_0)\big\||\nabla z|^2\big\|_{\Lo[6]}^2\nonumber\\\nonumber&\leq C_5\left(\intomega\big|\nabla|\nabla z|^2\big|^2\right)^\nfrac{5}{6}\left(\intomega|\nabla z|^2\right)^\nfrac{1}{3}+C_5\left(\intomega|\nabla z|^2\right)^2\\&\leq \frac{1}{2}\intomega\big|\nabla|\nabla z|^2\big|^2+C_6M^2\quad\text{for all }t>t_0.
\end{align}
The estimation of the remaining term on the right in \eqref{eq:cond_n-l2-combineq} is more involved. First, note that by \eqref{global_gnb_constants1} we have
\begin{align*}
\intomega|\nabla z|^6\leq K_2\left(\intomega\big|\nabla|\nabla z|^2\big|^2\right)\left(\intomega|\nabla z|^2\right)+K_2\left(\intomega|\nabla z|^4\right)\left(\intomega|\nabla z|^2\right)\quad\text{for all }t>t_0,
\end{align*}
where additionally by the Cauchy-Schwarz inequality
$\intomega|\nabla z|^4\leq\left(\intomega|\nabla z|^6\right)^\nfrac{1}{2}\left(\intomega|\nabla z|^2\right)^\nfrac{1}{2}
$ for all $t>t_0$, so that an application of Young's inequality combined with our assumption \eqref{eq:cond_n-l2_nlnn-req} implies that
\begin{align*}
\intomega|\nabla z|^6&\leq K_2\left(\intomega\big|\nabla|\nabla z|^2\big|^2\right)\left(\intomega|\nabla z|^2\right)+\eta\intomega|\nabla z|^6+\frac{K_2^2}{4\eta}\left(\intomega|\nabla z|^2\right)^3\\
&\leq K_2M\intomega\big|\nabla|\nabla z|^2\big|^2+\eta\intomega|\nabla z|^6+\frac{K_2^2M^3}{4\eta}\quad\text{for all }t>t_0
\end{align*}
and therefore
\begin{align}\label{eq:cond_n-l2-eq4}
(8+8\eta)\intomega|\nabla z|^6\leq\frac{8(1+\eta)K_2M}{1-\eta}\intomega\big|\nabla|\nabla z|^2\big|^2+\frac{2(1+\eta)K_2^2M^3}{(1-\eta)\eta}\quad\text{for all }t>t_0.
\end{align}
Collecting \eqref{eq:cond_n-l2-eq1}--\eqref{eq:cond_n-l2-eq4}, we infer from \eqref{eq:cond_n-l2-combineq} that for some $C_8>0$ we have
\begin{align}\label{eq:cond_n-l2-eq5}
\frac{\intd}{\intd t}\left\{\intomega n^2+\intomega|\nabla z|^4\right\}+C_7\intomega|\nabla n|^2+C_7\intomega\big|\nabla|\nabla z|^2\big|^2\leq C_8\quad\text{for all }t>t_0,
\end{align}
where $C_7:=\min\left\{\frac{1}{2},2-2\eta-\frac{8(1+\eta)K_2M}{1-\eta}\right\}$ is positive due to \eqref{eq:cond_n-l2_Mbound}. In order to conclude the desired bounds, we want to derive from the inequality above a differential inequality of the form $y'(t)+Cy^2(t)\leq C$, where $y(t):=\intomega n^2(\cdot,t)+\intomega|\nabla z(\cdot,t)|^4$ and $C>0$. To this end, we still need to estimate the terms without time derivatives, arising in \eqref{eq:cond_n-l2-eq5} on the left, from below. By making use of the \GNI, we firstly obtain upon use of the mass conservation and \eqref{eq:cond_n-l2_m-req} that
\begin{align*}
\left(\intomega n^2\right)^2&\leq C_9\left(\intomega |\nabla n|^2\right)\left(\intomega n\right)^2+C_9\left(\intomega n\right)^4\leq C_9m_0^2\intomega |\nabla n|^2+C_9m_0^4\quad\text{for all }t>t_0
\end{align*}
for some $C_9>0$, and secondly, relying on \eqref{eq:cond_n-l2_nlnn-req}, we find $C_{10}>0$ such that
\begin{align*}
\left(\intomega|\nabla z|^4\right)^2&\leq C_{10}\left(\intomega\big|\nabla|\nabla z|^2\big|^2\right)\left(\intomega|\nabla z|^2\right)^2+C_{10}\left(\intomega|\nabla z|^2\right)^4\\&\leq C_{10}M^2\intomega\big|\nabla|\nabla z|^2\big|^2+C_{10}M^4	\quad\text{for all }t>t_0.
\end{align*}
Thus, letting $C_{11}:=\max\{2C_9m_0^2,2C_{10}M^2\}$, we see that $y$ satisfies
\begin{align*}
y'(t)+C_{12}y^2(t)\leq C_{13}\quad\text{for all }t>t_0,
\end{align*}
with $C_{12}:=\frac{C_7}{C_{11}}$ and $C_{13}:=C_8+\frac{C_9m_0^4+C_{10}M^4}{C_{11}}$. By application of an ODE comparison argument, we observe that $\bar{y}(t):=\frac{2}{C_{12}(t-t_0)}+\sqrt{\frac{2C_{13}}{C_{12}}}$ satisfies $y(t)\leq \bar{y}(t)$ for all $t>t_0$, implying that 
\begin{align*}
y(t)\leq \frac{2}{C_{12}\tau}+\sqrt{\frac{2C_{13}}{C_{12}}}\quad\text{for all }t\geq t_0+\tau
\end{align*}
and thus proving \eqref{eq:cond_n-l2_n-l2}.
\end{bew}

\subsection{Eventual smoothness for generalized solutions with small mass}\label{sec4:evs}

For the our next proof we will require the following result demonstrated in \cite[Lemma 2.6]{Wang2016}, which is based on an application the Trudinger-Moser inequality combined with a spatio-temporal estimate on $\nabla\ln(n_\eps+1)$ in $L^2$.

\begin{lemma}\label{lem:moser_neps}
There exists $K_4>0$ such that for all $\eps\in(0,1)$ the solution to \eqref{CSz}--\,\eqref{ICz} satisfies
\begin{align*}
\intot\ln\bigg\{\frac{1}{|\Omega|}\intomega(n_\eps(x,s)+1)^2\intd x\bigg\}\intd s\leq K_4\left(1+\intomega n_0\right)t+K_4\left(\intomega z_0+\intomega n_0\right)\quad\text{for all }t>0.
\end{align*}
\end{lemma}

Relying on the properties previously established for $\F$, we can now determine some possibly large time $\ts$ depending on the initial data. But not on $\eps\in(0,1)$, for which $\intomega n_\eps|\ln n_\eps|$, $\intomega|\nabla z_\eps|^2$ and $\F(n_\eps,z_\eps)$ are sufficiently small for all times beyond $\ts$. This in turn will then ensure that we can obtain the conditional estimates featured in Section \ref{sec4:obtainingbound} for times larger than $\ts$.

\begin{lemma}\label{lem:neps-lnneps_nabzeps-l2-bound}
Let $K_2,K_3$ be as in \eqref{global_gnb_constants1} and \eqref{global_gnb_constants3}, respectively. There exist constants $\ms,\Gamma,M>0$ and $\mu\in(0,1)$ such that
\begin{align}\label{eq:neps-lnneps-GammaMreq}
\Gamma<\frac{1}{4K_3}-\frac{\mu|\Omega|}{e}\qquad\text{and}\qquad M<\frac{1}{4K_2},
\end{align}
and such that if the initial data $(n_0,c_0,u_0)$ satisfy \eqref{IR} as well as
\begin{align}\label{eq:neps-lnneps-massneps-bound}
m:=\intomega n_0\leq \ms,
\end{align}
then one can find $\ts>0$ such that for each $\eps\in(0,1)$ the solution $(n_\eps,z_\eps,u_\eps)$ of \eqref{CSz}--\,\eqref{ICz} satisfies
\begin{align}\label{eq:neps-lnneps-Fbound}
\F\big(n_\eps(\cdot,t),z_\eps(\cdot,t)\big)\leq \Gamma\quad\text{for all }t\geq\ts,
\end{align}
and
\begin{align}\label{eq:neps-lnneps-nbound}
\intomega n_\eps(\cdot,t)\left|\ln n_\eps(\cdot,t)\right|\leq\frac{1}{4K_3}+\frac{2|\Omega|}{e}\quad\text{for all }t\geq\ts,
\end{align}
as well as 
\begin{align}\label{eq:neps-lnneps-zbound}
\intomega|\nabla z_\eps(\cdot,t)|^2\leq M\quad\text{for all }t\geq\ts.
\end{align}
\end{lemma}

\begin{bew}
We fix $M\in\big(0,\frac{1}{4K_2}\big)$ and afterwards choose some small $\mu\in(0,1)$, such that 
\begin{align}\label{eq:neps-lnneps-Mbound}
\frac{2\mu|\Omega|}{e}\leq\frac{M}{2}\qquad\text{and}\qquad 0<\frac{1}{4K_3}-\frac{\mu|\Omega|}{e}.
\end{align}
Upon these choices, we can pick $\Gamma>0$ fulfilling the first inequality in \eqref{eq:neps-lnneps-GammaMreq} as well as
\begin{align}\label{eq:neps-lnneps-gammaM}
\Gamma\leq\frac{M}{4}.
\end{align}
Furthermore, letting $K_4$ be provided by Lemma \ref{lem:moser_neps} we can find $\eta\in(0,1)$ such that
\begin{align}\label{eq:omega-e-pwr-req}
\eta|\Omega|e^{16K_4}\leq\frac{\Gamma}{4}.
\end{align}
Relying on the previous choices and with $K_3,K_u$ given by \eqref{global_gnb_constants3} and Lemma \ref{lem:l4_bound_u}, respectively, we introduce the positive number
\begin{align}\label{eq:neps-lnneps-mstar-req}
\ms:=\min\left\{1,\frac{\Gamma}{4\ln\tfrac{1}{\eta\mu}},\frac{\Gamma}{8},\frac{1}{5K_3^2K_{u}|\Omega|^\frac{1}{4}}\right\},
\end{align}
where the positivity follows from the facts $\mu,\eta<1$. Now given $(n_0,c_0,u_0)$ such that \eqref{IR} and \eqref{eq:neps-lnneps-massneps-bound} hold, we find $\ell>0$ such that $\intomega|u_0|^4\leq\ell$, due to $\DA\hookrightarrow\Lo[4]$ (\cite[Lemma 2.3 iv)]{caolan16_smalldatasol3dnavstokes}). Moreover, since $\lambda_1>0$, we can easily find $t_0\geq0$ such that 
\begin{align}\label{eq:neps-lnneps-mstar-t0}
\ell e^{-\lambda_1 t_0}+\ms\leq\frac{1}{4K_3^2K_u|\Omega|^\frac{1}{4}}
\end{align}
holds. We next claim that the asserted inequalities are true if we fix some large $\ts$ satisfying the conditions
\begin{align}\label{eq:neps-lnneps-tstar-cond}
(1+m)\ts\geq \intomega z_0+m,\qquad m\ts\geq\intomega z_0,\qquad\text{and}\qquad \ts>2t_0,
\end{align}
with $z_0$ as defined in \eqref{ICz}. To verify this claim we define the sets
\begin{align*}
S_1(\eps):=\left\{t\in(0,\ts)\,\Big\vert\,\ln\Big\{\frac{1}{|\Omega|}\intomega(n_\eps(\cdot,t)+1)^2\Big\}>8K_4(1+m)\right\}
\end{align*} 
and
\begin{align*}
S_2(\eps):=\left\{t\in(0,\ts)\,\Big\vert\,\intomega|\nabla z_\eps(\cdot,t)|^2>8m\right\}
\end{align*}
and estimate their respective sizes. By Lemma \ref{lem:moser_neps} we know that for all $\eps\in(0,1)$ we have
\begin{align*}
I_1(\eps):=\int_0^{\ts}\ln\Big\{\frac{1}{|\Omega|}\intomega(n_\eps(\cdot,t)+1)^2\Big\}\intd t\leq K_4(1+m)\ts+K_4\left(\intomega z_0+m\right),
\end{align*}
so that the first condition in \eqref{eq:neps-lnneps-tstar-cond} combined with our definition of $S_1(\eps)$ shows that
\begin{align*}
2K_4(1+m)\ts&\geq K_4(1+m)\ts+K_4\left(\intomega z_0+m\right)\geq I_1(\eps)\geq 8K_4(1+m)|S_1(\eps)|
\end{align*}
holds for all $\eps\in(0,1)$, meaning that
\begin{align}\label{eq:neps-lnneps-S1-size}
|S_1(\eps)|\leq\frac{\ts}{4}\quad\text{for all }\eps\in(0,1).
\end{align}
In pursuance of a similar bound for the size of $|S_2(\eps)|$, we recall that by Lemma \ref{lem:bound_zeps} we have
\begin{align*}
I_2(\eps):=\int_0^{\ts}\!\!\intomega|\nabla z_\eps|^2\leq\intomega z_0+m\ts\quad\text{for all }\eps\in(0,1).
\end{align*}
Relying on the second inequality in \eqref{eq:neps-lnneps-tstar-cond} and the definition of $S_2(\eps)$ we infer that
\begin{align*}
2m\ts\geq\intomega z_0+m\ts\geq I_2(\eps)\geq 8m|S_2(\eps)|
\end{align*}
holds for all $\eps\in(0,1)$ and hence
\begin{align}\label{eq:neps-lnneps-S2-size}
|S_2(\eps)|\leq\frac{\ts}{4}\quad\text{for all }\eps\in(0,1).
\end{align}
Now, \eqref{eq:neps-lnneps-S1-size} and \eqref{eq:neps-lnneps-S2-size} guarantee that
\begin{align*}
\big|(0,\ts)\setminus\!\big(S_1(\eps)\cup S_2(\eps)\big)\big|\geq\frac{\ts}{2}\quad\text{for all }\eps\in(0,1),
\end{align*}
so that we conclude from the third inequality in \eqref{eq:neps-lnneps-tstar-cond} that for any $\eps\in(0,1)$ we can pick some $t_\eps\in(t_0,\ts)$ such that 
\begin{align}\label{eq:neps-lnneps-ineqs1}
\ln\left\{\frac{1}{|\Omega|}\intomega\big(n_\eps(\cdot,t_\eps)+1\big)^2\right\}\leq 8K_4(1+m)\qquad\text{and}\qquad\intomega|\nabla z_\eps(\cdot,t_\eps)|^2\leq 8m
\end{align}
hold. Relying on the elementary estimate $s\ln\frac{s}{\mu}\leq\eta(s+1)^2+s\ln\frac{1}{\eta\mu}$ for all $s>0$ (c.f. \cite[Lemma 5.5]{Win16CS2}), we can combine the mass conservation from Remark \ref{rem:masscons-sys-csz} with \eqref{eq:neps-lnneps-massneps-bound} and the first part of \eqref{eq:neps-lnneps-ineqs1} to obtain that
\begin{align*}
\intomega n_\eps(\cdot,t_\eps)\ln\frac{n_\eps(\cdot,t_\eps)}{\mu}\leq\eta\intomega\big(n_\eps(\cdot,t_\eps)+1\big)^2+\ln\frac{1}{\eta\mu}\intomega n_\eps(\cdot,t_\eps)\leq\eta|\Omega| e^{8K_4(1+m)}+m\ln\frac{1}{\eta\mu}.
\end{align*}
Now, recalling the first and second requirement for $\ms$ from \eqref{eq:neps-lnneps-mstar-req}, as well as \eqref{eq:omega-e-pwr-req}, we see that 
\begin{align*}
\intomega n_\eps(\cdot,t_\eps)\ln\frac{n_\eps(\cdot,t_\eps)}{\mu}\leq\eta|\Omega|e^{16K_4}+m\ln\frac{1}{\eta\mu}\leq\frac{\Gamma}{4}+\frac{\Gamma}{4}=\frac{\Gamma}{2}.
\end{align*}
In a similar fashion, the third part of \eqref{eq:neps-lnneps-ineqs1} in conjunction with the second inequality contained in \eqref{eq:neps-lnneps-mstar-req} entails that
\begin{align*}
\frac{1}{2}\intomega|\nabla z_\eps(\cdot,t_\eps)|^2\leq\frac{\Gamma}{2}
\end{align*}
and thus we obtain that
\begin{align*}
\F\big(n_\eps(\cdot,t_\eps),z_\eps(\cdot,t_\eps)\big)=\intomega n_\eps(\cdot,t_\eps)\ln\frac{n_\eps(\cdot,t_\eps)}{\mu}+\frac{1}{2}\intomega|\nabla z_\eps(\cdot,t_\eps)|^2\leq \Gamma.
\end{align*}
In accordance with \eqref{eq:neps-lnneps-GammaMreq} and \eqref{eq:neps-lnneps-mstar-t0}, this allows for the application of Lemma \ref{lem:decreasing-energy}, implying that
\begin{align}\label{eq:neps-lnneps-gammabound}
\F\big(n_\eps(\cdot,t),z_\eps(\cdot,t)\big)\leq\Gamma\quad\text{for all }t\geq t_\eps,
\end{align}
which, since $t_\eps<\ts$, immediately establishes \eqref{eq:neps-lnneps-Fbound} again due to \eqref{eq:neps-lnneps-GammaMreq}.
Now, to verify that also \eqref{eq:neps-lnneps-nbound} and \eqref{eq:neps-lnneps-zbound} hold, we recall that in view of Lemma \ref{lem:Fprop} we have
\begin{align*}
\intomega n_\eps(\cdot,t)|\ln n_\eps(\cdot,t)|\leq \F\big(n_\eps(\cdot,t),z_\eps(\cdot,t)\big)+\ln\mu\intomega n_\eps(\cdot,t)+\frac{2|\Omega|}{e}.
\end{align*}
Therefore, \eqref{eq:neps-lnneps-gammabound}, the fact $\mu<1$ and once more \eqref{eq:neps-lnneps-GammaMreq} imply 
\begin{align*}
\intomega n_\eps(\cdot,t)|\ln n_\eps(\cdot,t)|\leq\Gamma+\frac{2|\Omega|}{e}<\frac{1}{4K_3}+\frac{2|\Omega|}{e}\quad\text{for all }t\geq t_\eps,
\end{align*}
proving \eqref{eq:neps-lnneps-nbound}, because $\ts>t_\eps$. Similarly, again relying on Lemma \ref{lem:Fprop} and \eqref{eq:neps-lnneps-gammabound}, we conclude that due to \eqref{eq:neps-lnneps-gammaM} and the first restriction in \eqref{eq:neps-lnneps-Mbound}, we have
\begin{align*}
\intomega|\nabla z_\eps(\cdot,t)|^2\leq 2\F\big(n_\eps(\cdot,t),z_\eps(\cdot,t)\big)+\frac{2\mu|\Omega|}{e}\leq 2\Gamma+\frac{2\mu|\Omega|}{e}\leq\frac{M}{2}+\frac{M}{2}= M\quad\text{for all }t\geq t_\eps,
\end{align*}
which proves \eqref{eq:neps-lnneps-zbound}.
\end{bew}

The bounds for $\intomega n_\eps\ln n_\eps$ and $\intomega|\nabla z_\eps|^2$ at hand, we can first draw on the conditional estimates on $\intomega|\nabla z_\eps|^4$ from Section \ref{sec4:obtainingbound} and afterwards on the conditional regularity estimates from Section \ref{sec4:regularity} to obtain the following result.

\begin{proposition}\label{prop:evsmooth}
Let $\ms\!>0$ be as provided by Lemma \ref{lem:neps-lnneps_nabzeps-l2-bound}. Suppose that $(n_0,c_0,u_0)$ satisfy \eqref{IR} as well as
\begin{align*}
\intomega n_0\leq\ms,
\end{align*}
and let $(n,c,u)$ denote the global generalized solution of \eqref{CN}--\,\eqref{IC} from Theorem \ref{thm:globsol}. Then there exists $T>0$ such that
\begin{align}\label{eq:prop-evsmooth-evreg}
n\in\CSp{2,1}{\bomega\times[T,\infty)},\quad c\in\CSp{2,1}{\bomega\times[T,\infty)}\quad\text{and}\quad u\in\CSp{2,1}{\bomega\times[T,\infty);\R^2},
\end{align}
that
\begin{align*}
c(x,t)>0\quad\text{for all }x\in\bomega\text{ and any }t\geq T,
\end{align*}
and such that $(n,c,u)$ solves \eqref{CN}--\,\eqref{IC} classically in $\Omega\times(T,\infty)$. Moreover, one can find $\mu>0$ such that
\begin{align}\label{eq:prop-evfuncbound}
\F\big(n(\cdot,t),z(\cdot,t)\big)<\frac{1}{4K_3}-\frac{\mu|\Omega|}{e}\quad\text{for all }t\geq T,
\end{align}
with $z:=-\ln\frac{c}{\|c_0\|_{\Lo[\infty]}}$.
\end{proposition}

\begin{bew}
Let $K_2,K_3$ be provided by \eqref{global_gnb_constants1} and \eqref{global_gnb_constants3}, respectively. In view of Lemma \ref{lem:neps-lnneps_nabzeps-l2-bound} we can find $\mu\in(0,1)$, $\Gamma\in\big(0,\frac{1}{4K_3}-\frac{\mu|\Omega|}{e}\big)$, $M\in(0,\frac{1}{4K_2})$, $L>0$ and $\ts>0$ such that for any choice of $\eps\in(0,1)$ we have
\begin{align}\label{eq:prop-evsmooth-epsfuncbound}
\F\big(n_\eps(\cdot,t),z_\eps(\cdot,t)\big)&\leq \Gamma\quad\text{for all }t>\ts
\end{align}
and
\begin{align*}
\intomega n_\eps(\cdot,t)|\ln n_\eps(\cdot,t)|\leq L\quad\text{as well as}\quad\intomega|\nabla z_\eps(\cdot,t)|^2\leq M\quad\text{for all }t>\ts.
\end{align*}
Since $M<\frac{1}{4K_2}$, we may employ Lemma \ref{lem:cond_n-l2_nabz-l4_bound} to obtain $C_1>0$ such that for any $\eps\in(0,1)$ we have
\begin{align*}
\intomega|\nabla z_\eps(\cdot,t)|^4\leq C_1\quad\text{for all }t>\ts+1.
\end{align*}
This bound at hand, Lemma \ref{lem:cond-reg-hoelder2} yields $\theta\in(0,1)$ such that for each $T>\ts+2$ we can pick $C_2(T)>0$ such that
\begin{align*}
\|n_\eps\|_{\CSp{2+\theta,1+\frac{\theta}{2}}{\bomega\times[\ts+2,T]}}+\|z_\eps\|_{\CSp{2+\theta,1+\frac{\theta}{2}}{\bomega\times[\ts+2,T]}}+\|u_\eps\|_{\CSp{2+\theta,1+\frac{\theta}{2}}{\bomega\times[\ts+2,T]}}\leq C_2(T)
\end{align*}
for all $\eps\in(0,1)$. In view of the Arzelà-Ascoli theorem, we can find a subsequence $(\eps_{j_k})_{k\in\N}$ of the sequence provided by Lemma \ref{lem:convergence1}, along which $n_\eps$, $z_\eps$ and $u_\eps$ are convergent in $\CSploc{2,1}{\bomega\times[\ts+2,\infty)}$. The respective limits of $n_\eps$, $z_\eps$ and $u_\eps$ must clearly coincide with $n$, $z$ and $u$, which ensures that $n$, $c$ and $u$ have the desired regularity properties in \eqref{eq:prop-evsmooth-evreg}. Additionally, the continuity of $z$ implies $c>0$ in $\bomega\times[T,\infty)$ and passing to the limit for $\eps=\eps_{j_k}\searrow0$ in \eqref{eq:prop-evsmooth-epsfuncbound} we easily obtain \eqref{eq:prop-evfuncbound} due to $\Gamma<\frac{1}{4K_3}-\frac{\mu|\Omega|}{e}$. Letting $\eps=\eps_{j_k}\searrow0$ in \eqref{CSz} we first conclude that $(n,z,u)$ solves \eqref{CSznoeps}--\eqref{CSznoepsBC} with $f(\xi)\equiv\xi$ classically in $\Omega\times(T,\infty)$, which then in combination with $c>0$ in $\bomega\times[T,\infty)$ entails that $(n,c,u)$ solve \eqref{CN}--\eqref{IC} classically in $\Omega\times[T,\infty)$.
\end{bew}

\subsection{Stabilization of solutions with small energy}\label{sec4:stabilization}
This section discusses the last missing part for the proof of Theorem \ref{thm:evsmooth}, which is the convergence properties featured therein. Since from the last section we already known, that our generalized solutions will be classical solutions after some waiting time, we will concern our investigation only with convergence of classical solutions to \eqref{CSznoeps}. Before proving the desired large time behavior we require one additional preparation in form of a time-independent Hölder bound from $\nabla z$.

\begin{lemma}\label{lem:time-indep-hoelder-nabzeps}
For all $m_0>0$, $M>0$, $\tau>0$ there exist $\theta\in(0,1)$ and $C>0$ such that if for $f\in C^{3}\!\left([0,\infty)\right)$ satisfying \eqref{eq:fnoepsprop} and $t_0\geq0$ the triple $(n,z,u)\in\CSp{0}{\bomega\times[t_0,\infty)}\cap\CSp{2,1}{\bomega\times(t_0,\infty)}$ is a classical solution of \eqref{CSznoeps}--\,\eqref{CSznoepsBC} in $\Omega\times(t_0,\infty)$  satisfying
\begin{align*}
\intomega n(\cdot,t_0)\leq m_0,
\end{align*}
and
\begin{align*}
\intomega|\nabla z(\cdot,t)|^4\leq M\quad\text{for all }t>t_0,
\end{align*}
it holds that
\begin{align}\label{eq:timidephoeld-zhoel}
\|\nabla z(\cdot,t)\|_{\CSph{\theta}{\bomega}}\leq C\quad\text{for all }t\geq t_0+\tau.
\end{align}
\end{lemma}

\begin{bew}
The arguments are quite similar to the ones employed in \cite[Lemma 4.9]{Win16CS2} and we will not recount all details here. First, we note that by Lemma \ref{lem:cond-reg-n-inf} we can find $C_1>0$ such that 
\begin{align}\label{eq:timidephoel-nlp}
\|n(\cdot,t)\|_{\Lo[4]}\leq C_1\quad\text{for all }t\geq \overline{t_0}:=t_0+\frac{\tau}{2}.
\end{align}
Now, we may choose some $\beta\in(0,1)$ close to 1 such that $\beta>\frac{1}{4}$ and afterwards $q>1$ satisfying $\frac{1}{4}<\frac{1}{q}<\frac{5}{4}-\beta$. With these values fixed we will make use of several well knwon estimates for the Neumann heat semigroup $\big(e^{-s B}\big)_{s\geq0}$ in $\Lo[4]$, where $B:=-\Delta+1$ (e.g. \cite{win10jde}). Moreover, for any fixed $\theta\in(0,2\beta-\frac{3}{2})$ we have that $D\!\left(B^\beta\right)\hookrightarrow \CSp{1+\theta}{\bomega}$ (\cite[Theorem 1.6.1]{hen81}) and hence
\begin{align}\label{eq:timidephoel-emb}
\|\nabla\varphi\|_{\CSph{\theta}{\bomega}}\leq C_2\|B^\beta\varphi\|_{\Lo[4]}\quad\text{for all }\varphi\in\D\!\left(B^\beta\right),
\end{align}
with some $C_2>0$. Letting 
\begin{align*}
S_1:=\max_{t\in[\overline{t_0},\overline{t_0}+1]}(t-\overline{t_0})^\beta\|\nabla z(\cdot,t)\|_{\CSph{\theta}{\bomega}}\quad\text{and}\quad S_2(T):=\max_{t\in[\overline{t_0}+1,T]}\|\nabla z(\cdot,t)\|_{\CSph{\theta}{\bomega}}
\end{align*}
for $T>\overline{t_0}+1$ we continue by estimating $S(T):=\max\left\{S_1,S_2(T)\right\}$. Consequently, with $t_1(t):=\max\{t-1,\overline{t_0}\}$ we start by representing $z(\cdot,t)$ according to
\begin{align}\label{eq:timidephoel-zrep}
z(\cdot,t)&=\overline{z(\cdot,t_1)}+e^{t-t_1}e^{-(t-t_1)B}\Big(z(\cdot,t_1)-\overline{z(\cdot,t_1)}\Big)-\int_{t_1}^t e^{t-s}e^{-(t-s)B}|\nabla z(\cdot,s)|^2\intd s\nonumber\\&\quad+\int_{t_1}^te^{t-s}e^{-(t-s)B}f\big(n(\cdot,s)\big)\intd s-\int_{t_1}^te^{t-s}e^{-(t-s)B}u(\cdot,s)\nabla z(\cdot,s)\intd s.
\end{align}
In the case of $t-\overline{t_0}\leq 1$ we make use of Young's inequality, \eqref{eq:timidephoel-emb}, the semigroup estimates for the Neumann heat semigroup, and the fact that $f(s)\leq s$ for all $s\geq0$ to obtain $C_3>0$ such that
\begin{align}\label{eq:timidephoel-zsemiest}
\nonumber\|\nabla z(\cdot,t)\|_{\CSph{\theta}{\bomega}}&\leq C_3e(t-\overline{t_0})^{-\beta}\|z(\cdot,\overline{t_0})-\overline{z(\cdot,\overline{t_0})}\|_{\Lo[4]}+C_3e\int_{\overline{t_0}}^t(t-s)^{-\gamma}\big\||\nabla z(\cdot,s)|^2\big\|_{\Lo[q]}\intd s\\&+C_3e\int_{\overline{t_0}}^t(t-s)^{-\beta}\|n(\cdot,s)\|_{\Lo[4]}\intd s+C_3e\int_{\overline{t_0}}^t(t-s)^{-\beta}\|u(\cdot,s)\|_{\Lo[2p]}^2\intd s,	
\end{align}
holds for all $t\leq\bar{t_0}+1$, where $\gamma:=\beta+\frac{1}{q}-\frac{1}{4}<1$. Herein, \eqref{eq:timidephoel-nlp} and Lemma \ref{lem:stokes_est}, and the fact that $\beta<1$ imply the existence of $C_4>0$ such that
\begin{align*}
C_3e\int_{\overline{t_0}}^t(t-s)^{-\beta}\|n(\cdot,s)\|_{\Lo[4]}\intd s+C_3e\int_{\overline{t_0}}^t(t-s)^{-\beta}\|u(\cdot,s)\|_{\Lo[8]}^2\intd s\leq C_4\int_{\overline{t_0}}^t(t-s)^{-\beta}\intd s\leq \frac{C_4}{1-\beta},
\end{align*}
for all $t\geq \bar{t_0}+1$, and the Poincaré inequality provides $C_5>0$ satisfying
\begin{align*}
\|z(\cdot,s)-\overline{z(\cdot,s)}\|_{\Lo[4]}\leq C_5\|\nabla z(\cdot,s)\|_{\Lo[4]}\leq C_5 M^\frac{1}{4}\quad\text{for all }s\geq\overline{t_0}.
\end{align*}
Furthermore, by means of the Hölder inequality we see that
\begin{align*}
\big\||\nabla z(\cdot,s)|^2\big\|_{\Lo[q]}\leq\|\nabla z(\cdot,s)\|_{\Lo[4]}^\frac{4}{q}\|\nabla z(\cdot,s)\|_{\Lo[\infty]}^a\leq M^\frac{1}{q}\|\nabla z(\cdot,s)\|_{\CSph{\theta}{\bomega}}^a\quad\text{for all }s\geq \overline{t_0},
\end{align*}
with $a:=\frac{2q-4}{q}$, and hence for all $t\geq\bar{t_0}+1$ we have
\begin{align*}
\int_{\overline{t_0}}^t\!(t-s)^{-\gamma}\big\||\nabla z(\cdot,s)|^2\big\|_{\Lo[q]}\!\intd s\leq M^\frac{1}{q} S_1^a(t-\overline{t_0})^{1-\gamma-\beta a}\!\!\int_0^1\!\!(1-\sigma)^{-\gamma}\sigma^{-\beta a}\intd\sigma
\leq C_6M^\frac{1}{q}S_1^a(t-\overline{t_0})^{1-\gamma-\beta a},
\end{align*}
where we used that $
\int_0^1(1-\sigma)^{-\gamma}\sigma^{-\beta a}\intd\sigma=:C_6$ is finite due to the facts that $0<a<1$, $0<\beta<1$ and $\gamma<1$.
Accordingly, from \eqref{eq:timidephoel-zsemiest} we infer that
\begin{align*}
(t-\overline{t_0})^{\beta}\|\nabla z(\cdot,t)\|_{\CSph{\theta}{\bomega}}\leq C_3C_5eM^\frac{1}{4}+C_3C_6eM^\frac{1}{q}S_1^a(t-\overline{t_0})^{1-\gamma+(1-a)\beta}+\frac{C_4}{1-\beta}\leq C_7+C_7S_1^a
\end{align*}
for all $t\in[\overline{t_0},\overline{t_0}+1]$, with some $C_7>0$, which implies that $S_1\leq\max\{1,(2C_7)^{\frac{1}{1-a}}\}$. Similarly, in the case $t\in[\overline{t_0},T]$ we conclude from \eqref{eq:timidephoel-zrep} that
\begin{align*}
\|\nabla z(\cdot,t)\|_{\CSph{\theta}{\bomega}}\leq C_8M^\frac{1}{4}+C_8M^\frac{1}{q}\int_{t-1}^t(t-s)^{-\gamma}\|\nabla z(\cdot,s)\|_{\CSph{\theta}{\bomega}}^a\intd s+C_8\int_{t-1}^t(t-s)^{-\beta}\intd s,
\end{align*}
for some $C_8>0$. In both of the cases $t\leq\overline{t_0}+2$ and $t>\overline{t_0}+2$ we can estimate
\begin{align*}
\int_{t-1}^t(t-s)^{-\gamma}\|\nabla z(\cdot,s)\|_{\CSph{\theta}{\bomega}}^a\intd s&\leq S_1^a\int_{t-1}^t(t-s)^{-\gamma}(s-\overline{t_0})^{-\beta a}\intd s+S_2^a(T)\int_{t-1}^t(t-s)^{-\gamma}\intd s\\
&\leq C_5S_1^a+\frac{1}{1-\gamma}S_2^a(T)
\end{align*}
with $C_5$ as defined above. Therefore, for suitable large $C_9>0$ we have
\begin{align*}
S_2(T)\leq C_9+C_9 S_2^a(T)\quad\text{for all }T>\overline{t_0}+1,
\end{align*}
which implies that $S_2(T)\leq\max\{1,(2C_9)^\frac{1}{1-a}\}=:S_2$ for all $T>\overline{t_0}+1$. Consequently, together with the previous estimate for $S_1$, this establishes \eqref{eq:timidephoeld-zhoel} with $C:=\max\{S_1,\frac{S_1}{\tau},S_2\}$.
\end{bew}

Assuming that the energy $\F(n,z)$ remains small for all times succeeding some waiting $T\geq0$, which according to Proposition \ref{prop:evsmooth} is true for the generalized solutions with small mass, we will now show that any given solution to \eqref{CSznoeps}--\eqref{CSznoepsBC} in $\Omega\times(T,\infty)$ will satisfy the asymptotic properties described in Theorem \ref{thm:evsmooth}. Here we explicitly allow $T=0$, because if the energy is already suitably small initially we can transfer these asymptotic properties also to the global classical solutions discussed in Section \ref{sec4:globclass}.

\begin{proposition}\label{prop:stabilization}
Assume $T\geq0$, $\ell>0$ and let $\ms>0$ be as in Lemma \ref{lem:neps-lnneps_nabzeps-l2-bound}. Suppose that for $f\in C^{3}\!\left([0,\infty)\right)$ satisfying \eqref{eq:fnoepsprop} the triple $(n,z,u)\in\CSp{0}{\bomega\times[T,\infty)}\cap\CSp{2,1}{\bomega\times(T,\infty)}$ is a classical solution of \eqref{CSznoeps}--\,\eqref{CSznoepsBC} in $\Omega\times(T,\infty)$ satisfying $z\in\CSp{0}{[T,\infty);\W[1,2]}$, $m:=\intomega n(\cdot,T)<\ms$, $0\leq n\not\equiv0$, and $\intomega|u(\cdot,T)|^4\leq \ell$, as well as
\begin{align}\label{eq:stabil-infF-prop}
\inf_{t>T}\F\big(n(\cdot,t),z(\cdot,t)\big)<\frac{1}{4K_3}-\frac{\mu|\Omega|}{e}
\end{align}
for some $\mu>0$. Then
\begin{align}\label{eq:stabil-nconv}
n(\cdot,t)\to \overline{n_{T}}:=\frac{1}{|\Omega|}\intomega n(\cdot,T)\quad\text{in }\Lo[\infty]\quad\text{as }t\to\infty,
\end{align}
and
\begin{align}\label{eq:stabil-nabzconv}
\nabla z(\cdot,t)\to 0\quad\text{in }\Lo[\infty]\quad\text{as }t\to\infty,
\end{align}
and
\begin{align}\label{eq:stabil-zconv}
\inf_{x\in\Omega}z(x,t)\to \infty\quad\quad\text{as }t\to\infty,
\end{align}
as well as
\begin{align}\label{eq:stabil-uconv}
u(\cdot,t)\to 0\quad\text{in }\Lo[\infty]\quad\text{as }t\to\infty.
\end{align}
\end{proposition}

\begin{bew}
The convergence of $n$ and $z$ can be proved by relying on the methods shown in \cite[Lemma 6.1]{Win16CS2}, whereas the decay of $u$ then follows by adapting the arguments illustrated in \cite[Lemma 5.3]{win_ct_fluid_3d}. For the sake of completeness we only recount the main steps and refer to the mentioned sources for more details. 
Recalling that $\ms<\tfrac{1}{4K_3^2K_{u}|\Omega|^\frac{1}{4}}$, we can first find $t_0>T$ such that
$\ell e^{-\lambda_1 (t_0-T)}+\ms\leq\frac{1}{4K_3^2K_{u}|\Omega|^\frac{1}{4}}$
and then rely on \eqref{eq:stabil-infF-prop} and Lemma \ref{lem:decreasing-energy} to see that we can pick $\ts>t_0>T$ such that
\begin{align}\label{eq:prop-stabi-eq05}
\frac{\intd}{\intd t}\F\big(n(\cdot,t),z(\cdot,t)\big)\leq 0\quad\text{for all }t>\ts,
\end{align}
and
\begin{align}\label{eq:prop-stabi-eq1}
\F\big(n(\cdot,t),z(\cdot,t)\big)<C_1:=\frac{1}{4K_3}-\frac{\mu|\Omega|}{e}\quad\text{for all }t>\ts,
\end{align}
and that with some $\kappa>0$,
\begin{align}\label{eq:prop-stabi-eq1.5}
\int_{\ts}^\infty\!\intomega\frac{|\nabla n|^2}{n}+\kappa\int_{\ts}^\infty\!\intomega|\Delta z|^2\leq C_2:=\frac{1}{4K_3}.
\end{align}
Since $(n,z,u)$ solve \eqref{CSznoeps} classically in $\Omega\times(T,\infty)$ by Remark \ref{rem:masscons-sys-csz} we have
\begin{align}\label{eq:prop-stabi-masscons}
\intomega n(\cdot,t)=m\quad\text{for all }t>T,
\end{align}
and thus, making use of \eqref{eq:Fineq_n} and \eqref{eq:prop-stabi-eq1}, we see that
\begin{align}\label{eq:prop-stabi-nlnbound}
\intomega n(\cdot,t)|\ln n(\cdot,t)|\leq \F\big(n(\cdot,t),z(\cdot,t)\big)+\ln\mu\intomega n(\cdot,t)+\frac{2|\Omega|}{e}\leq C_1+m\ln\mu+\frac{2|\Omega|}{e}
\end{align}
holds for all $t>\ts$. Since $\W[1,1]\hookrightarrow\Lo[2]$, a Poincaré--Sobolev inequality implies the existence of $C_3>0$ such that
\begin{align}\label{eq:prop-stabi-eq2}
\|\varphi-\overline{\varphi}\|_{\Lo[2]}\leq C_3\|\nabla \varphi\|_{\Lo[1]}\quad\text{for all }\varphi\in\W[1,1].
\end{align}
Similarly, by means of elliptic regularity theory we can find $C_4>0$ satisfying
\begin{align}\label{eq:prop-stabi-eq3}
\|\nabla \varphi\|_{\Lo[2]}\leq C_4\|\Delta\varphi\|_{\Lo[2]}\quad\text{for all }\varphi\in\W[2,2]\text{ such that }\frac{\partial\varphi}{\partial\nu}=0\text{ on }\romega.
\end{align}
According to \eqref{eq:prop-stabi-eq2} and the Cauchy-Schwarz inequality we thus have
\begin{align*}
\int_{\ts}^\infty\|n(\cdot,t)-\overline{n_T}\|^2_{\Lo[2]}\intd t\leq C_3^2\int_{\ts}^\infty\|\nabla u\|_{L[1]}^2\intd t\leq m C_3^2\int_{\ts}^\infty\!\intomega\frac{|\nabla n|^2}{n},
\end{align*}
whereas \eqref{eq:prop-stabi-eq3} shows that
\begin{align*}
\int_{\ts}^T\|\nabla z(\cdot,t)\|_{\Lo[2]}^2\intd t\leq C_4^2\int_{\ts}^\infty\!\intomega|\Delta z|^2.
\end{align*}
By combination of the two previous estimates with \eqref{eq:prop-stabi-eq1.5} we thereby see that
\begin{align}\label{eq:prop-stabi-eq4}
\int_{\ts}^\infty\left\{\|n(\cdot,t)-\overline{n_T}\|_{\Lo[2]}^2+\|\nabla z(\cdot,t)\|_{\Lo[2]}^2\right\}\intd t\leq C_2\big(mC_3^2+\frac{C_4^2}{\kappa}\big)
\end{align}
which implies that there must exist $(t_k)_{k\in\N}\subset(\ts,\infty)$ such that $t_k\to\infty$ and such that
\begin{align}\label{eq:prop-stabi-nzl2}
n(\cdot,t_k)\to\overline{n_T}\quad\text{in }\Lo[2]\quad\text{and}\quad \nabla z(\cdot,t_k)\to0\quad\text{in }\Lo[2]
\end{align}
as $k\to\infty$. Relying on the convexity of $0<\xi\mapsto \xi\ln \xi$ and the Jensen inequality we see that
\begin{align*}
\intomega\varphi\ln\varphi\intd x\geq\intomega\overline{\varphi}\ln\overline{\varphi}\quad\text{for all positive }\varphi\in\CSp{0}{\bomega},
\end{align*}
and thus, we can make use of the mean value theorem, the Cauchy-Schwarz inequality, the first convergence in \eqref{eq:prop-stabi-masscons}, and \eqref{eq:prop-stabi-nzl2} to obtain
\begin{align}\label{eq:prop-stabi-mvest}
0&\leq \intomega n(\cdot,t_k)\ln n(\cdot,t_k)-\intomega\overline{n_T}\ln\overline{n_T}=\intomega n(\cdot,t_k)\big(\ln n(\cdot,t_k)-\ln\overline{n_T}\big)\nonumber\\
&\leq\int_{\{n(\cdot,t_k)>\overline{n_T}\}}n(\cdot,t_k)\big(\ln n(\cdot,t_k)-\ln\overline{n_T}\big)\nonumber\\
&\leq\frac{1}{\overline{n_T}}\|n(\cdot,t_k)\|_{\Lo[2]}\|n(\cdot,t_k)-\overline{n_T}\|_{\Lo[2]}\to 0\quad\text{as }k\to\infty.
\end{align}
This, together with the definition of $\F$ and the second convergence established in \eqref{eq:prop-stabi-nzl2} shows that
$\F\big(n(\cdot,t_k),z(\cdot,t_k)\big)\to C_5:=\intomega \overline{n_T}\ln\frac{\overline{n_T}}{\mu}$ as $k\to\infty$, which in turn by the monotonicity property \eqref{eq:prop-stabi-eq05} implies
\begin{align*}
\F\big(n(\cdot,t),z(\cdot,t)\big)\to C_5\quad\text{as }t\to\infty.
\end{align*}
In view of \eqref{eq:prop-stabi-mvest} this convergence actually yields
\begin{align}\label{eq:prop-stabi-eq45}
\limsup_{t\to\infty}\intomega|\nabla z(\cdot,t)|^2=2\limsup_{t\to\infty}\left\{\F\big(n(\cdot,t),z(\cdot,t)\big)-\intomega n(\cdot,t)\ln\frac{n(\cdot,t)}{\mu}\right\}\leq 2C_5-2C_5=0.
\end{align}
Combining this with the bound provided by \eqref{eq:prop-stabi-nlnbound} we may first employ Lemma \ref{lem:cond_n-l2_nabz-l4_bound} and afterwards Lemma \ref{lem:cond-reg-nu-hoelder} and Lemma \ref{lem:time-indep-hoelder-nabzeps} to obtain $\tss>\ts$, $\theta\in(0,1)$ and $C_6>0$ such that
\begin{align}\label{eq:prop-stabi-eq5}
\|n\|_{\CSp{\theta,\frac{\theta}{2}}{\bomega\times[t,t+1]}}\leq C_6,\quad \|u\|_{\CSp{\theta,\frac{\theta}{2}}{\bomega\times[t,t+1]}}\leq C_6,\quad\text{and}\quad\|\nabla z(\cdot,t)\|_{\CSph{\theta}{\bomega}}\leq C_6
\end{align}
for all $t\geq\tss$. If the asserted convergence for $n$ in \eqref{eq:stabil-nconv} was false we could find $(\tilde{t}_k)_{k\in\N}\subset(\tss,\infty)$ and $C_7>0$ such that $\tilde{t}_k\to\infty$ as $k\to\infty$ and 
\begin{align*}
\|n(\cdot,\tilde{t}_k)-\overline{n_T}\|_{\Lo[\infty]}\geq C_7\quad\text{for all }k\in\N,
\end{align*}
implying that, due to the uniform convergence of $n$ in $\bomega\times[\tss,\infty)$ asserted by \eqref{eq:prop-stabi-eq5}, there exist $(x_k)_{k\in\N}\subset\Omega$, $r>0$, and $\tau>0$ such that $B_r(x_k)\subset\Omega$ for all $k\in\N$ and
\begin{align*}
\big|n(x,t)-\overline{n_T}\big|\geq\frac{C_7}{2}\quad\text{for all }x\in B_r(x_k)\text{ and each }t\in(\tilde{t}_k,\tilde{t}_k+\tau). 
\end{align*}
In turn this would show that
\begin{align*}
\int_{\tilde{t}_k}^{\tilde{t}_k+\tau}\|n(\cdot,t)-\overline{n_T}\|_{\Lo[2]}^2\intd t\geq \tau\frac{C_7^2}{4}\pi r^2\quad\text{for all }k\in\N, 
\end{align*}
contradicting the spatial-temporal estimate \eqref{eq:prop-stabi-eq4} and thus proving \eqref{eq:stabil-nconv}. In a similar fashion, assuming that \eqref{eq:stabil-nabzconv} is false, in view of the second portion of \eqref{eq:prop-stabi-eq5}, we could find $(\hat{t}_k)_{k\in\N}\subset(\tss,\infty)$, $(\hat{x}_k)_{k\in\N}\subset\Omega$, $r>0$, and $C_8>0$ such that $\hat{t}_k\to\infty$ as $k\to\infty$ and $B_r(\hat{x}_k)\subset\Omega$ for all $k\in\N$ as well as
\begin{align*}
|\nabla z(x,\hat{t}_k)|\geq C_8\quad\text{for all }x\in B_r(\hat{x}_k)\text{ and each }k\in\N.
\end{align*}
This implies that
\begin{align*}
\intomega|\nabla z(\cdot,\hat{t}_k)|^2\geq C_8^2\pi r^2\quad\text{for all }k\in\N,
\end{align*}
which contradicts \eqref{eq:prop-stabi-eq45} and thereby proves \eqref{eq:stabil-nabzconv}. For \eqref{eq:stabil-zconv} we make use of the fact that \eqref{eq:stabil-nconv} together with the nontriviality of $n$ establishes the existence of some $t_{\star\star\star}>T$ satisfying
\begin{align*}
n(x,t)>\frac{\overline{n_T}}{2}\quad\text{for all }x\in\Omega\text{ and }t>t_{\star\star\star},
\end{align*}
whence, by relying on the nonnegativity of $z$ and parabolic comparison with the function $\bomega\times[t_{\star\star\star},\infty)\ni(x,t)\mapsto\tfrac{\overline{n_t}}{2}(t-t_{\star\star\star})$, we see that
\begin{align*}
z(x,t)\geq\frac{\overline{n_T}}{2}(t-t_{\star\star\star})\quad\text{for all }x\in\Omega\text{ and }t>t_{\star\star\star},
\end{align*}
ensuring \eqref{eq:stabil-zconv}. In order to prove \eqref{eq:stabil-uconv}, we recall that the Stokes operator $A$ in $L^2_\sigma\!\left(\Omega\right)$ is positive and self-adjoint with compact inverse and as such, there exists a complete orthonormal basis $(\psi_k)_{k\in\N}$ of eigenfunctions of $A$ to positive eigenvalues $\lambda_k$, $k\in\N$. Since $\bigcup_{m\in\N}\text{span}\left\{\psi_k|k\leq m\right\}$ is dense in $L^2_\sigma\!\left(\Omega\right)$, in view of the uniform Hölder continuity of $u$ in $\Omega\times(\tss,\infty)$ from \eqref{eq:prop-stabi-eq5}, we only have to show that for each $k\in\N$ we have 
\begin{align}\label{eq:prop-stabi-eq55}
\intomega u(x,t)\cdot\psi_k(x)\intd x\to0\quad\text{as }t\to\infty.
\end{align}
To this end we fix $k\in\N$ and let $y(t):=\intomega u(x,t)\cdot\psi_k(x)\intd x$, $t>T$. From the third equation in \eqref{CSznoeps}, the eigenfunction property of $\psi_k$, as well as the fact that $\dive\psi_k=0$ we obtain
\begin{align}\label{eq:prop-stabi-eq6}
y'(t)=-\lambda_k\intomega u\cdot\psi_k+\intomega \big(n-\overline{n_T}\big)\nabla\phi\cdot\psi_k\quad\text{for all }t>T.
\end{align}
Since $n\to\overline{n_T}$ in $\Lo[\infty]$ as $t\to\infty$ by \eqref{eq:stabil-nconv}, for any given $\delta>0$ we can find $t_{\diamond}>T$ such that 
\begin{align*}
\left|\intomega\big(n(x,t)-\overline{n_T}\big)\nabla\phi\cdot\psi_k(x)\intd x\right|\leq \frac{\delta\lambda_k}{2}\quad\text{for all }t>t_{\diamond},
\end{align*}
which shows upon integration of \eqref{eq:prop-stabi-eq6} that, due to the boundedness of $u$ in $\Omega\times(T,\infty)$, we have
\begin{align*}
y(t)&<y(t_{\diamond})e^{-\lambda_k(t-t_{\diamond})}+\frac{\lambda_k \delta}{2}\int_{t_{\diamond}}^te^{-\lambda_k(t-s)}
<C_9e^{-\lambda_k(t-t_{\diamond})}+\frac{\delta}{2}\quad\text{for all }t>t_{\diamond},
\end{align*}
with some $C_9>0$. Now letting $t_{\diamond\diamond}:=\max\left\{t_{\diamond},t_{\diamond}+\frac{1}{\lambda_k}\ln\frac{2C_9}{\delta}\right\}$ we have
\begin{align*}
|y(t)|<\delta\quad\text{for all }t>t_{\diamond\diamond},
\end{align*}
yielding \eqref{eq:prop-stabi-eq55} and thus completing the proof.
\end{bew}

All that is left is to gather the results of our previous two propositions to conclude the proof of Theorem \ref{thm:evsmooth}.

\begin{proof}[\textbf{Proof of Theorem \ref{thm:evsmooth}:}]
With $\ms>0$ provided by Lemma \ref{lem:neps-lnneps_nabzeps-l2-bound} we obtain from Proposition \ref{prop:evsmooth} that for any initial data $(n_0,c_0,u_0)$ satisfying \eqref{IR} as well as \eqref{eq:critmass}, there exists $T>0$ such that the solution $(n,c,u)$ from Theorem \ref{thm:globsol} has the regularity properties featured in \eqref{eq:evreg} and the positivity of $c$ in $\bomega\times(T,\infty)$ as claimed in \eqref{eq:large-time-positivity-c} are valid. Since \eqref{eq:prop-evfuncbound} from Proposition \ref{prop:evsmooth} furthermore guarantees that $\inf_{t>T}\F\big(n(\cdot,t),z(\cdot,t)\big)<\tfrac{1}{4K_3}-\tfrac{\mu|\Omega|}{e}$, we may employ Proposition \ref{prop:stabilization} to obtain \eqref{eq:conv-n} and \eqref{eq:conv-nabc}.
\end{proof}

\subsection{Global classical solutions for small initial data. Proof of Theorem \ref{thm:smalldataglobalclass}}\label{sec4:globclass}
As mentioned in the introduction, the result featured in Theorem \ref{thm:smalldataglobalclass} is a by-product of our previous analysis. Our main tools in the proof will on one hand be the fact that the assumed smallness conditions for the initial data, expressed in \eqref{eq:init-small} and \eqref{eq:init-energ}, allows for the choice of $t_0=0$ in Lemma \ref{lem:decreasing-energy}, and on the other hand the uniqueness statement from Lemma \ref{lem:loc_ex_approx_c}. The uniqueness statement is essential, since we can only guarantee the global existence for our approximate solutions when $f(s)\equiv f_\eps(s)$ with $f_\eps(s)$ provided by \eqref{eq:feps_def}.

\begin{proof}[\textbf{Proof of Theorem \ref{thm:smalldataglobalclass}:}]
We denote by $(n,c,u)$ the local classical solution from Lemma \ref{lem:loc_ex_approx_c} for $f(s)\equiv s$, extended to its maximal existence time $\Tm\in(0,\infty]$. Then, writing $z:=-\ln\big(\frac{c}{\|c_0\|_{\Lo[\infty]}}\big)$ and $\tau:=\min\{1,\frac{\Tm}{2}\}$, we infer that $C_1:=\|n\|_{\LSp{\infty}{\Omega\times(0,\tau)}}$ is finite, by the continuity of $n$ in $\bomega\times[0,\Tm)$. On the other hand, let us also consider the approximate problems \eqref{CSz} and denote the corresponding solutions by $(n_\eps,z_\eps,u_\eps)$ with $\eps\in(0,1)$. According to \cite[Section 2.1]{Wang2016} these solutions are global for each of these $\eps\in(0,1)$. For these solutions and $\mu$ as in \eqref{eq:init-energ} we have
\begin{align*}
\F\big(n_\eps(\cdot,0),z_\eps(\cdot,0)\big)=C_2:=\intomega n_0\ln\frac{n_0}{\mu}+\frac{1}{2}\intomega\frac{|\nabla z_0|^2}{c_0^2}\quad\text{for all }\eps\in(0,1),
\end{align*}
and furthermore, defining $\mss:=\frac{1}{8K_3^2K_{u}|\Omega|^\frac{1}{4}}$ we conclude that the inequalities contained in \eqref{eq:init-small} imply
\begin{align*}
\intomega|u_0|^4 e^{-\lambda_1 t}+\intomega n_0<\frac{1}{4K_3^2K_{u}|\Omega|^\frac{1}{4}}\quad\text{for all }t>0.
\end{align*}
In light of \eqref{eq:mass_cons_n} and \eqref{eq:init-energ} we have $C_2<\frac{1}{4K_3}-\frac{\mu|\Omega|}{e}$, Lemma \ref{lem:decreasing-energy} becomes applicable, asserting that
\begin{align*}
\F\big(n_\eps(\cdot,t),z_\eps(\cdot,t)\big)\leq C_2\quad\text{for all }t>0\text{ and each }\eps\in(0,1).
\end{align*}
Thanks to Lemma \ref{lem:Fprop} this implies that for any $\eps\in(0,1)$ we have
\begin{align*}
\intomega n_\eps\left|\ln n_\eps(\cdot,t)\right|\leq C_2+\ln\mu\intomega n_0+\frac{2|\Omega|}{e}
\quad\text{and}\quad
\intomega |\nabla z_\eps|^2\leq M:=2C_2+\frac{2\mu|\Omega|}{e}\quad\text{for all }t>0.
\end{align*}
Herein, the second restriction on $C_2$ from \eqref{eq:init-energ} shows that
\begin{align*}
M<\frac{2}{8K_2}-\frac{2\mu|\Omega|}{e}+\frac{2\mu|\Omega|}{e}=\frac{1}{4K_2}.
\end{align*}
Hence, we may employ Lemma \ref{lem:cond_n-l2_nabz-l4_bound} to find $C_3>0$ such that
\begin{align*}
\intomega|\nabla z_\eps(\cdot,t)|^4\leq C_3\quad\text{for all }t>\frac{\tau}{2}\text{ and each }\eps\in(0,1).
\end{align*}
In turn, Lemma \ref{lem:cond-reg-n-inf} becomes applicable and provides $C_4>0$ such that
\begin{align}\label{eq:theo2-ninfbound}
\|n_\eps(\cdot,t)\|_{\Lo[\infty]}\leq C_4\quad\text{for all }t>\tau\text{ and every }\eps\in(0,1).
\end{align}
Now, fixing $\eps\in(0,1)$ so small such that it satisfies $\eps\leq\min\left\{\frac{1}{C_1},\frac{1}{C_4}\right\}$,
we see that by the definition of $f_\eps$ in \eqref{eq:feps_def} we have
\begin{align*}
f_\eps(n)=n\quad\text{in }\bomega\times[0,\tau],
\end{align*}
from which , in view of the uniqueness statement contained in Lemma \ref{lem:loc_ex_approx_c} when applied to the system \eqref{CSc} with $f\equiv f_\eps$, we infer that\vspace*{-5pt}
\begin{align*}
(n,z,u)\equiv(n_\eps,z_\eps,u_\eps)\quad\text{in }\bomega\times[0,\tau]\vspace*{-5pt}
\end{align*}
for our fixed $\eps$. On the other hand, relying on \eqref{eq:theo2-ninfbound} and the second restriction on $\eps$ we also have $f_\eps(n_\eps)\equiv n_\eps$ in $\bomega\times(\tau,\infty)$ and $(n_\eps,z_\eps,u_\eps)$ actually solves \eqref{CSznoeps} in $\Omega\times(\tau,\infty)$ with $f(s)\equiv s$. Now, making use of the uniqueness result from Lemma 2.1 once more, when applied to \eqref{CSc} with $f(s)\equiv s$, guarantees that $\Tm=\infty$ and that $(n,z,u)\equiv(n_\eps,z_\eps,u_\eps)$ in $\Omega\times(0,\infty)$. The desired convergence properties easily follow from Proposition \ref{prop:stabilization}, since $C_2<\frac{1}{4K_3}-\frac{\mu|\Omega|}{e}$.
\end{proof}
\vspace*{-16pt}
\section*{Acknowledgements}
\vspace*{-6pt}
The author acknowledges the support of the {\em Deutsche Forschungsgemeinschaft} in the context of the project
  {\em Analysis of chemotactic cross-diffusion in complex frameworks}. 

\vspace*{-6pt}
\footnotesize{

}

\begin{thebibliography}{40}
\providecommand{\natexlab}[1]{#1}
\providecommand{\url}[1]{\texttt{#1}}
\expandafter\ifx\csname urlstyle\endcsname\relax
  \providecommand{\doi}[1]{doi: #1}\else
  \providecommand{\doi}{doi: \begingroup \urlstyle{rm}\Url}\fi

\bibitem[Adler(1966)]{Adler708}
J.~Adler.
\newblock Chemotaxis in bacteria.
\newblock \emph{Science}, 153\penalty0 (3737):\penalty0 708--716, 1966.

\bibitem[Amann(2000)]{Amann00}
H.~Amann.
\newblock Compact embeddings of vector-valued {S}obolev and {B}esov spaces.
\newblock \emph{Glas. Mat. Ser. III}, 35(55)\penalty0 (1):\penalty0 161--177,
  2000.

\bibitem[Biler et~al.(1994)Biler, Hebisch, and Nadzieja]{Bil94}
P.~Biler, W.~Hebisch, and T.~Nadzieja.
\newblock The {D}ebye system: existence and large time behavior of solutions.
\newblock \emph{Nonlinear Anal.}, 23\penalty0 (9):\penalty0 1189--1209, 1994.

\bibitem[Cao and Lankeit(2016)]{caolan16_smalldatasol3dnavstokes}
X.~Cao and J.~Lankeit.
\newblock Global classical small-data solutions for a three-dimensional
  chemotaxis {N}avier-{S}tokes system involving matrix-valued sensitivities.
\newblock \emph{Calc. Var. Partial Differential Equations}, 55\penalty0
  (4):\penalty0 Paper No. 107, 39, 2016.

\bibitem[Evans(2010)]{evans}
L.~C. Evans.
\newblock \emph{Partial differential equations}, volume~19 of \emph{Graduate
  Studies in Mathematics}.
\newblock American Mathematical Society, Providence, RI, second edition, 2010.

\bibitem[Fujie et~al.(2014)Fujie, Ito, and Yokota]{FIY14}
K.~Fujie, A.~Ito, and T.~Yokota.
\newblock Existence and uniqueness of local classical solutions to modified
  tumor invasion models of {C}haplain-{A}nderson type.
\newblock \emph{Adv. Math. Sci. Appl.}, 24\penalty0 (1):\penalty0 67--84, 2014.

\bibitem[Fujie et~al.(2016)Fujie, Ito, Winkler, and Yokota]{FIYW16}
K.~Fujie, A.~Ito, M.~Winkler, and T.~Yokota.
\newblock Stabilization in a chemotaxis model for tumor invasion.
\newblock \emph{Discrete Contin. Dyn. Syst.}, 36\penalty0 (1):\penalty0
  151--169, 2016.

\bibitem[Giga and Sohr(1991)]{GigSohr91}
Y.~Giga and H.~Sohr.
\newblock Abstract {$L^p$} estimates for the {C}auchy problem with applications
  to the {N}avier-{S}tokes equations in exterior domains.
\newblock \emph{J. Funct. Anal.}, 102\penalty0 (1):\penalty0 72--94, 1991.

\bibitem[Henry(1981)]{hen81}
D.~Henry.
\newblock \emph{Geometric Theory of Semilinear Parabolic Equations}, volume 840
  of \emph{Lecture Notes in Mathematics}.
\newblock Springer Berlin Heidelberg, 1981.

\bibitem[Hillen and Painter(2009)]{HP09}
T.~Hillen and K.~J. Painter.
\newblock A user’s guide to {PDE} models for chemotaxis.
\newblock \emph{J. Math. Biol.}, 58\penalty0 (1-2):\penalty0 183--217, 2009.

\bibitem[Horstmann and Winkler(2005)]{HoWin05_bvblowchemo}
D.~Horstmann and M.~Winkler.
\newblock Boundedness vs. blow-up in a chemotaxis system.
\newblock \emph{J. Differential Equations}, 215\penalty0 (1):\penalty0 52 --
  107, 2005.

\bibitem[Keller and Segel(1971)]{KS71travbands}
E.~F. Keller and L.~A. Segel.
\newblock Traveling bands of chemotactic bacteria: A theoretical analysis.
\newblock \emph{J. Theor. Biol.}, 30\penalty0 (2):\penalty0 235 -- 248, 1971.

\bibitem[Kozono et~al.(2016)Kozono, Miura, and Sugiyama]{kozono15}
H.~Kozono, M.~Miura, and Y.~Sugiyama.
\newblock Existence and uniqueness theorem on mild solutions to the
  {K}eller-{S}egel system coupled with the {N}avier-{S}tokes fluid.
\newblock \emph{J. Funct. Anal.}, 270\penalty0 (5):\penalty0 1663--1683, 2016.

\bibitem[Lady\v{z}enskaja et~al.(1968)Lady\v{z}enskaja, Solonnikov, and
  Ural'ceva]{LSU}
O.~A. Lady\v{z}enskaja, V.~A. Solonnikov, and N.~N. Ural'ceva.
\newblock \emph{Linear and quasilinear equations of parabolic type}.
\newblock Translations of mathematical monographs. American Mathematical
  Society, 1968.

\bibitem[Lankeit(2016)]{Lan16_M3AS}
J.~Lankeit.
\newblock Long-term behaviour in a chemotaxis-fluid system with logistic
  source.
\newblock \emph{Math. Models Methods Appl. Sci.}, 26\penalty0 (11):\penalty0
  2071--2109, 2016.

\bibitem[Li and Zhao(2015)]{li2015initial}
H.~Li and K.~Zhao.
\newblock Initial--boundary value problems for a system of hyperbolic balance
  laws arising from chemotaxis.
\newblock \emph{J. Differential Equations}, 258\penalty0 (2):\penalty0
  302--338, 2015.

\bibitem[Li et~al.(2014)Li, Li, and Wang]{LiLiWang14}
J.~Li, T.~Li, and Z.-A. Wang.
\newblock Stability of traveling waves of the {K}eller-{S}egel system with
  logarithmic sensitivity.
\newblock \emph{Math. Models Methods Appl. Sci.}, 24\penalty0 (14):\penalty0
  2819--2849, 2014.

\bibitem[Lorz(2010)]{lorz10}
A.~Lorz.
\newblock Coupled chemotaxis fluid model.
\newblock \emph{Math. Mod. Meth. Appl. S.}, 20\penalty0 (06):\penalty0
  987--1004, 2010.

\bibitem[Mizoguchi and Souplet(2014)]{MS14}
N.~Mizoguchi and P.~Souplet.
\newblock Nondegeneracy of blow-up points for the parabolic {K}eller-{S}egel
  system.
\newblock \emph{Ann. Inst. H. Poincar\'e Anal. Non Lin\'eaire}, 31\penalty0
  (4):\penalty0 851--875, 2014.

\bibitem[Nagai and Ikeda(1991)]{NagIke91}
T.~Nagai and T.~Ikeda.
\newblock Traveling waves in a chemotactic model.
\newblock \emph{J. Math. Biol.}, 30\penalty0 (2):\penalty0 169--184, 1991.

\bibitem[Porzio and Vespri(1993)]{PorzVesp93}
M.~M. Porzio and V.~Vespri.
\newblock H\"older estimates for local solutions of some doubly nonlinear
  degenerate parabolic equations.
\newblock \emph{J. Differential Equations}, 103\penalty0 (1):\penalty0
  146--178, 1993.

\bibitem[Quittner and Souplet(2007)]{QS07}
P.~Quittner and P.~Souplet.
\newblock \emph{Superlinear parabolic problems}.
\newblock Birkh\"auser Advanced Texts: Basler Lehrb\"ucher. Birkh\"auser
  Verlag, Basel, 2007.

\bibitem[Rosen(1978)]{ROSEN1978}
G.~Rosen.
\newblock Steady-state distribution of bacteria chemotactic toward oxygen.
\newblock \emph{Bull. Math. Biol.}, 40\penalty0 (5):\penalty0 671 -- 674, 1978.

\bibitem[Sell and You(2002)]{sellyou}
G.~R. Sell and Y.~You.
\newblock \emph{Dynamics of evolutionary equations}, volume 143 of
  \emph{Applied Mathematical Sciences}.
\newblock Springer-Verlag, New York, 2002.

\bibitem[Sohr(2001)]{sohr}
H.~Sohr.
\newblock \emph{The {N}avier-{S}tokes equations}.
\newblock Birkh\"auser Advanced Texts: Basler Lehrb\"ucher. Birkh\"auser
  Verlag, Basel, 2001.

\bibitem[Tao(2011)]{Tao-consumption_JMAA11}
Y.~Tao.
\newblock Boundedness in a chemotaxis model with oxygen consumption by
  bacteria.
\newblock \emph{J. Math. Anal. Appl.}, 381\penalty0 (2):\penalty0 521--529,
  2011.

\bibitem[Tao and Winkler(2012)]{TaoWin-evsmooth_JDE12}
Y.~Tao and M.~Winkler.
\newblock Eventual smoothness and stabilization of large-data solutions in a
  three-dimensional chemotaxis system with consumption of chemoattractant.
\newblock \emph{J. Differential Equations}, 252\penalty0 (3):\penalty0
  2520--2543, 2012.

\bibitem[Tao et~al.(2013)Tao, Wang, and Wang]{TWWDCDSB13}
Y.~Tao, L.~Wang, and Z.-A. Wang.
\newblock Large-time behavior of a parabolic-parabolic chemotaxis model with
  logarithmic sensitivity in one dimension.
\newblock \emph{Discrete Contin. Dyn. Syst. Ser. B}, 18:\penalty0 821--845,
  2013.

\bibitem[Tuval et~al.(2005)Tuval, Cisneros, Dombrowski, Wolgemuth, Kessler, and
  Goldstein]{tuval2005bacterial}
I.~Tuval, L.~Cisneros, C.~Dombrowski, C.~W. Wolgemuth, J.~O. Kessler, and R.~E.
  Goldstein.
\newblock Bacterial swimming and oxygen transport near contact lines.
\newblock \emph{Proc. Natl. Acad. Sci. U.S.A.}, 102\penalty0 (7):\penalty0
  2277--2282, 2005.

\bibitem[Wang(2016)]{Wang2016}
Y.~Wang.
\newblock Global large-data generalized solutions in a two-dimensional
  chemotaxis-{S}tokes system with singular sensitivity.
\newblock \emph{Boundary Value Problems}, 2016\penalty0 (1):\penalty0 177,
  2016.

\bibitem[Wang and Xiang(2015)]{Wang20157578}
Y.~Wang and Z.~Xiang.
\newblock Global existence and boundedness in a {K}eller–{S}egel–{S}tokes
  system involving a tensor-valued sensitivity with saturation.
\newblock \emph{J. Differential Equations}, 259\penalty0 (12):\penalty0 7578 --
  7609, 2015.

\bibitem[Wang(2013)]{Wang13surv}
Z.-A. Wang.
\newblock Mathematics of traveling waves in chemotaxis---review paper.
\newblock \emph{Discrete Contin. Dyn. Syst. Ser. B}, 18\penalty0 (3):\penalty0
  601--641, 2013.

\bibitem[Wang et~al.(2016)Wang, Xiang, and Yu]{Wang20162225}
Z.-A. Wang, Z.~Xiang, and P.~Yu.
\newblock Asymptotic dynamics on a singular chemotaxis system modeling onset of
  tumor angiogenesis.
\newblock \emph{J. Differential Equations}, 260\penalty0 (3):\penalty0 2225 --
  2258, 2016.

\bibitem[Winkler(2010)]{win10jde}
M.~Winkler.
\newblock Aggregation vs. global diffusive behavior in the higher-dimensional
  {K}eller-{S}egel model.
\newblock \emph{J. Differential Equations}, 248\penalty0 (12):\penalty0
  2889--2905, 2010.

\bibitem[Winkler(2012)]{win_fluid_final}
M.~Winkler.
\newblock Global large-data solutions in a chemotaxis-({N}avier-){S}tokes
  system modeling cellular swimming in fluid drops.
\newblock \emph{Comm. Partial Differential Equations}, 37\penalty0
  (2):\penalty0 319--351, 2012.

\bibitem[Winkler(2015{\natexlab{a}})]{win15_chemorot}
M.~Winkler.
\newblock Large-data global generalized solutions in a chemotaxis system with
  tensor-valued sensitivities.
\newblock \emph{SIAM J. Math. Anal.}, 47\penalty0 (4):\penalty0 3092--3115,
  2015{\natexlab{a}}.

\bibitem[Winkler(2015{\natexlab{b}})]{win_ct_fluid_3d}
M.~Winkler.
\newblock Boundedness and large time behavior in a three-dimensional
  chemotaxis-{S}tokes system with nonlinear diffusion and general sensitivity.
\newblock \emph{Calc. Var. Partial Differential Equations}, 54\penalty0
  (4):\penalty0 3789--3828, 2015{\natexlab{b}}.

\bibitem[Winkler(2016{\natexlab{a}})]{Win16CS1}
M.~Winkler.
\newblock The two-dimensional {K}eller-{S}egel system with singular sensitivity
  and signal absorption: global large-data solutions and their relaxation
  properties.
\newblock \emph{Math. Models Methods Appl. Sci.}, 26\penalty0 (5):\penalty0
  987--1024, 2016{\natexlab{a}}.

\bibitem[Winkler(2016{\natexlab{b}})]{Win16CS2}
M.~Winkler.
\newblock The two-dimensional {K}eller-{S}egel system with singular sensitivity
  and signal absorption: {E}ventual smoothness and equilibration of small-mass
  solutions.
\newblock 2016{\natexlab{b}}.
\newblock Preprint.

\bibitem[Winkler(2017)]{win15_chemonavstokesfinal}
M.~Winkler.
\newblock How far do chemotaxis-driven forces influence regularity in the
  {N}avier-{S}tokes system?
\newblock \emph{Trans. Amer. Math. Soc.}, 369\penalty0 (5):\penalty0
  3067--3125, 2017.

\end{thebibliography}
\end{document}